\documentclass{cocv}

\usepackage[normalem]{ulem}
\usepackage{lettrine}
\usepackage{boondox-cal}
\usepackage[utf8]{inputenc}
\usepackage{amsmath}
\usepackage{amsfonts}
\usepackage{amssymb}
\usepackage{showkeys}
\usepackage{amsthm}
\usepackage{bm}
\usepackage{xcolor}
\usepackage[english]{babel}
\usepackage[T1]{fontenc}
\usepackage{mathrsfs}
\usepackage{graphicx}
\graphicspath{{./Images/}}
\usepackage{setspace}
\usepackage{bigints}
\newcommand{\rn}{\mathbb{R}^n}
\newcommand{\essinf}{\text{essinf}}

\newcommand{\ol}{\overline}
\newcommand{\cl}[1]{\mathcal{#1}}

\newcommand{\extpbf}[1]{ \mathcal{#1}_{\bm{\varepsilon}} }
\newcommand{\opprocess}{(\overline{S},\overline{w}^0,\overline{w},\overline{\alpha},\overline{y}^0,\overline{y},\overline{\beta})}

\newcommand{\dist}[2]{\text{dist}\left(#1,#2\right)}

\DeclareMathOperator*{\Diff}{Diff}
\DeclareMathOperator*{\Lin}{Lin}

\def\prj{\mathfrak{pr}}

\def\rr{\mathbb{R}}
\definecolor{mydarkred}{RGB}{180, 10, 10}

\usepackage{cancel}
\def\bel{\begin{equation}\label}
	\def\eeq{\end{equation}}
\def\ds{\displaystyle}

\def\setmap{\rightsquigarrow}

\numberwithin{equation}{section}

\begin{document}
%%-----------------------------
%%      the top matter
%%-----------------------------
\title{Goh conditions 
for minima of\\ nonsmooth problems with unbounded controls.}
%\thanks{EVENTUALI RINGRAZIAMENTI}

\author{Francesca Angrisani}\address{Department of Mathematics "T.Levi-Civita", University of Padova, Italy}
\author{Franco Rampazzo}\sameaddress{1}
\date{...}
\begin{abstract} Higher order necessary conditions for a minimizer of an optimal control problem  are generally obtained for systems whose dynamics is continuously differentiable in the state variable.  Here, by making use of the notion of set-valued Lie bracket we obtain a Goh-type condition for a control affine system with Lipschitz continuous dynamics and unbounded controls. In order to manage the simultaneous lack of smoothness of the adjoint equation and of the Lie bracket-like variations we make use of the notion of Quasi Differential Quotient.{ We conclude the paper with a worked out example
	where }the established higher order condition is capable to rule out the optimality of a control verifying the standard maximum principle.  \end{abstract}
\begin{resume} 
	Les conditions nécessaires d'ordre supérieur pour un minimiseur d'un problème de contrôle optimal sont généralement obtenues pour des systèmes dont la dynamique est $C^1$ en la variable d'état. Ici, en utilisant la notion de crochet de Lie multivoque, nous obtenons une condition de type Goh pour un système de contrôle affine avec une dynamique Lipschitz et des contrôles non bornés. Afin de gérer le manque de régularité simultané de l'équation adjointe et des variations de type crochet de Lie, nous utilisons la notion de Quasi Differential Quotient.  Nous concluons le papier avec un exemple qui montre comment la condition d'ordre supérieur établie permet d'exclure l'optimalité d'un contrôle vérifiant le principe du maximum classique. \end{resume}
\maketitle

\section{Introduction}
The so-called higher order necessary conditions for minima,  a classical subject of investigation since the early developments of Calculus of Variations, have been variously generalized to Optimal Control Theory. In particular, they constitute a crucial issue in geometric control, in connection with the Lie algebraic structure associated with the dynamics.  Clearly, such a structure is present  only under hypotheses of $C^\infty$ regularity for the involved vector fields.  Yet, 
if only   Lie brackets up to a certain  length are used, one can assume  less demanding smoothness hypotheses, as it happens, for instance, for the Goh condition.  Indeed, the mere definition of  Lie bracket of two vector fields $g_i,g_j$, namely $[g_i,g_j]:=Dg_jg_i-Dg_ig_j$,  only requires that they be differentiable. 
It is then plausible to wonder if one can  further weaken  the regularity assumptions, for instance by allowing vector fields to be just locally Lipschitz continuous. On the one hand, this would be  of obvious interest for applications, and, on the other hand, it would be in line with the vast and  rich literature  on {\it nonsmooth} optimal control, which, since the early Seventies, has generated some {\it maximum principles} involving suitable   notions of generalized differentiation (see \cite{AubinFrankowska},\cite{Clarke}, \cite{Mordukhovich} \cite{Sussmann2},\cite{SussmannWarga},\cite{Vinter},\cite{Warga1}, and \cite{Warga2}).

Let us observe that, when two vector fields  $g_i,g_j$  are  locally Lipschitz continuous, Rademacher's theorem implies that the differentiation domains $\Diff(g_i)$ and $\Diff(g_j)$ of  $g_i$ and $g_j$, respectively, have full measure, so that  the Lie bracket $[g_i,g_j]$ (existing on $\Diff(g_i) \cap \Diff(g_j)$) turns out to be defined almost everywhere. Actually, some results including Lie brackets  have been generalized to nonsmooth systems  by means of such an  almost everywhere approach (see, for example, \cite{RampSuss2007},\cite{Simic}, \cite{Montanari}, \cite {Cardin}).
As an alternative, and  in connection with {controllability} problems for Lipschitz continuous control systems,  an  {\it  everywhere-defined, set-valued,} Lie bracket  has been introduced in \cite{RampSuss2001}. If  one uses $\overline{\text{co}} (E)$ to  denote the closed convex hull of a given set $E$, this bracket is defined   by setting, for every $x\in\rr^n$,
\bel{lieset}
[g_i,g_j]_{set}(x)=\overline{\text{co}}\Big\{\lim\limits_{n\to \infty} [g_i,g_j](x_n),\,\,\quad \lim_{n\to \infty} x_n  = x,\,\, x_n \in \Diff(f)\cap \Diff(g)\Big\},
\eeq
where we mean that limits are taken along all sequences $(x_n)\subset \Diff(g_i)\cap \Diff(g_j)$  converging to $x$.
In \cite{Frobenius},  \cite{RampSuss2007} this notion of   set-valued Lie bracket --which can be defined on differential manifolds as well--  proved suitable for the generalization of  some basic results of differential geometry like the  commutativity criterion, Frobenius Theorem, and Rashevski-Chow Theorem.
Therefore, one might naturally consider the following question:
\vskip0.5truecm
\begin{itemize}
	\item[\bf Q.]
	{ \it Given a minimizer  of a nonsmooth optimal control problem, can one  complement the (first order) Maximum Principle with  conditions    which involve set-valued Lie brackets?} 
\end{itemize}
\vskip0.5truecm

In the present  paper we provide a first positive answer to question {\bf Q} in relation with  a  control affine  minimum problem
\bel{minpr}(P)\quad\begin{cases}
	{\text{minimize }} \Psi(T,x(T))\\
	\text{over processes $(T,u,x)$\,\, of} 
	\\\ds
	\frac{dx}{dt}=f(x)+\sum\limits_{i=1}^m g_i(x)u^i\qquad\\
	%	\displaystyle% \frac{d\nu}{dt}(t)=|u(t)|,\\
	x(0) =\hat{x}\quad
	%\quad \nu(T)\le K,
	\|u\|_1\leq K\quad
	(T,x(T)) \in \mathfrak{T}.
\end{cases}
\eeq
and  proving a necessary condition generalizing Goh condition. 
We assume that the vector fields \linebreak {$g_i:\rr^n \to \rr^n, \, i=1,\ldots,m$,  are locally  Lipschitz  continuous}, while  
the controls $u$ take values in $\rr^m$ and have   $L^1$-norm less than  or equal to $K$.   The 
end-time $T$ and the final state $x(T)$ are subject to a constraint of the form $(T,x(T))\in \mathfrak{T}$, where the {\it  target} $ \mathfrak{T}$ is a given closed subset of the time-space product $\rr_+\times\rr^n$. This is actually a simplified version of the problem addressed in the paper,  in which the drift $f$  is allowed to depend on   a bounded control $a$ as well  and, moreover, a Lagrangian  cost  $l=l(t,x,u)$ is considered together with the final cost $\Psi$. Moreover, the set where the  control $u$ takes values is not necessarily the whole $\rr^m$ but it might be allowed to be a   closed cone $C=C_1\times C_2$, where $C_2$ is a cone and for some non negative integer  $m_1\leq m$,  $C_1 \subseteq \xR^{m_1}$ is a closed cone containing the coordinate axes.
%For the sake of clarity, let us assume  that   the cost $\Psi$  is smooth.  {\angr  Adattare questa farse a seconda della scelta che faremo}.
%if we use $AC$ to denote the set of absolutely continuous paths (each one on some interval $[0,T]$), a control trajectory pair $((a,u),x) \in L^1\times AC$ can hardly be a minimizer of problem 
%\eqref{minpr}. Indeed, the latter is a slow growth problem, so that minimizer are often {\it impulsive.} In other words, 

Since the velocities of the trajectories are unbounded and  no suitable growth condition prevents the system from an impulsive behaviour, we first embed the problem in the time-space optimal control problem 
$$ (P_{ext})\quad\begin{cases}
	\text{minimize } \Psi(y^0(S),y(S)) \\
	\text{ over processes } (S,w^0,w,y^0,y),\\\hbox{where  $(y^0,y):[0,S]\to\rr^{1+n}$ solves}\\ \ds\frac{dy^0}{ds}=w^0
	\\ \displaystyle \frac{dy}{ds}=f(y)w^0+\sum\limits_{i=1}^m g_i(y)w^i\qquad
	\\ \displaystyle (y^0,y)(0)=(0,\hat{x})\quad \|w\|_1\leq K, \quad (t,y)(S)\in  \mathfrak{T},\end{cases} $$ where the  controls $(w^0,w)$ belong to the set 
$\bigcup\limits_{S>0} \Big\{(w^0,w) \in L^\infty([0,S],\rr_+\times \rr^m):  w^0(s)+|w(s)|= 1 \Big\}$ and $y^0$ stands for the actual time $t$ variable.
Problem $(P_{ext})$ is simply obtained from $(P)$ by first reparametrizing  time through $$\displaystyle t(s)= y^0(s) :=\int_0^s w^0(\sigma)d\sigma, \qquad w^0>0, \qquad w(s):=(u\circ t)(s) w^0(s), \qquad y(s):=(x\circ y^0)(s),$$ and then allowing  subintervals  $I\subseteq[0,S]$ such that  $w^0(s)\equiv 0 \,\forall s\in I$ ({\it impulsive subintervals}).   Notice, in particular, that, unlike $(P)$,   $(P_{ext})$ is a problem with $L^\infty$-bounded controls. Incidentally, notice also that, since $(P_{ext})$ is rate-independent,  the control constraint $ w^0+|w|= 1$
is not restrictive.

As for necessary conditions for a minimizer  $(\overline{S},\overline{w}^0,\overline{w},\overline{y}^0,\overline{y})$ of the extended problem $(P)_{ext}$, answering question {\bf Q}  should mean complementing the usual, non smooth,   maximum principle (in one of the available versions)    with  conditions  that tell something about the relation between the corresponding adjoint variable $p(\cdot)$ and set-valued Lie brackets $[g_i,g_j]_{set}$.

Our main result, which we  
state below in a simplified form --see Theorem \ref {TeoremaPrincipale} for a rigorous and more general statement--, actually says that at almost every $s\in [0,\ol{S}]$ the following non smooth Goh condition holds true: 

\begin{thrm}[\bf Maximum Principle]\label{TeoremaPrincipaleintro}
	Let $(\overline{S},\overline{w}^0,\overline{w},\overline{y}^0,
	\overline{y})$ be a  local minimizer for the extended  problem 
	$(P_{ext})$, and assume that $\|\ol{w}\|_1<K$. 
	
	Then there exist multipliers $(p_0,p,\lambda) \in \rr^*\times AC([0,\ol{S}];(\rr^n)^*) \times \mathbb{R^*}$ such that, besides  the standard necessary conditions of Pontryagin Maximum Principle (expressed withing a suitable non-smooth setting), we have
	\bel{condII}
	\boxed{0\in p(s) \cdot  [g_i,g_j]_{set}(\ol{y}(s)),}
	\eeq
	for any $i,j \in \{1,\ldots,m\}$ and for almost any $s \in [0,\ol{S}]$.
\end{thrm}
Clearly, the importance of such a result relies on the possibility that   a given  control $\hat u$ while being allowed by the standard, first order, maximum principle, does not verify condition \eqref{condII}, which would rule out the optimality of $\hat u$.  A toy example at the end of the paper illustrates this circumstance. 

Let us mention that a crucial  tool for the proof of Theorem \ref{TeoremaPrincipaleintro} (in the more general version of Theorem \ref{TeoremaPrincipale}), is represented by the notion of Quasi Differential Quotient ($QDQ$), a  generalized differentation (valid  for set-valued maps) introduced in \cite{RampPal} as a special case of H. Sussmann's Approximate Generalized Differential Quotients \cite{Sussmann}.
Actually, this tool and the corresponding notion of approximating multicone are flexible enough to allow the managing, whitin the same theoretical frame, of two different  kinds of nonsmoothness: the one generated by the adjoint inclusion (which involves Clarke's generalized Jacobian), and the one derived by the utilization of set-valued brackets. In particular, one  
constructs variational $QDQ$-approximating multicones generated  by both multiple 
set-valued Lie brackets and the classical needle variations, and proves a linear separability condition with $QDQ$-approximating multicones of the set of profitable states. 
In turn, this geometrical  result is equivalent to the existence of the multipliers $(p_0,p,\lambda) $ stated in  Theorem \ref{TeoremaPrincipaleintro}.

To conclude this presentation, let us briefly mention that the interest for the issues treated in this paper is justified by several applications, for instance  in classical mechanics (as soon one identifies   the  control with a moving part of a given mechanical system),  in neurological dynamics, or aerospace navigation \cite{BressanAldo1989,BressanRampazzo2010,GRR,Catlla2008,Marle,AzimovBishop2005, Rampazzo1991,BressanMotta1993}. A further application which  deserves some attention is  that of driftless control systems, in particular the case of sub-Riemannian geometry. For instance,  a result as the one presented here might be regarded as a contribution  to the study of  sub-Riemannian metrics having low regularity (see, for instance, \cite{Montanari}).
%Finally, we refer the reader to the last section for some comments concerning possible developments in the  direction of considering  the Goh condition with bounded controls  as well as  other Lie-bracket involving higher order conditions, still under low-regularity assumptions.
\vskip0.5truecm

 Finally, let us spend a few words on possible developments of the present work. On the one hand,  aiming to  higher order necessary conditions involving {\it iterated set-valued  Lie-Brackets}  ( which were introduced in \cite{RampFeleqi}), one should begin by proving  that  such brackets are $QDQ$ of suitable multiflows. This is not straighforward, for iterated  set-valued  Lie-Brackets are larger then the objects one would obtain by a mere recursive approach. 
On the other hand, it might be interesting  to generalize the Goh-like condition proved here  to the standard  case of {\it bounded} controls. Actually, in this case  the construction of the bracket-like variation cannot ignore the  influence of the non-zero drift (as instead it happens with unbounded controls). A similar, though probably more complex, argument  would  obviously be involved in the investigation of a non-smooth version of the Legendre-Clebsch condition.

\vskip0.8truecm
The paper is organized as follows: in the next subsection we introduce some conventions for notation; Section 2 is devoted to the introduction of some   crucial  theoretical  tools, like set-valued Lie brackets   and Quasi Differential Quotients;  in Section 3 we  present the  minimum problem as well as its  extended version, and state the main result (Theorem \ref{TeoremaPrincipale});  Section 4 is dedicated to the proof of the main result, while in Section 5 we describe a toy example where the effectiveness of the new, higher order, necessary condition is displayed in ruling out a first order extremal.
%  The paper ends with a short subsection containing some ideas for future work.
\subsection{Notation} \label{notations}
{
	We shall use $\xR_+$ and $\xR_-$ to denote {the intervals} $[0,+\infty[$ and $]-\infty,0],$ respectively. The space of all linear operators from a vector space $X$ to a vector space $Y$ will be denoted by $\Lin(X,Y)$. The elements of an Euclidean space $\xR^q$, $q\geq 0$, will be thought as {column vectors}, while {\ row vectors} will {stand for} the linear one-forms, i.e. the elements of the dual space {$(\xR^q)^*:=\Lin(\xR^q,\xR)$}. If $p,q$ are positive integers, the space  $\Lin(\xR^p,\xR^q)$  will be sometimes regarded as the  space of {$q\times p$ real-valued} matrices.}
For every $i=1,\ldots,q$,  $\mathbf{e}_i$ [resp. $\mathbf{e}^i$]  will denote the $i$-th vector of the canonical basis of $\xR^q$ [resp.  $(\xR^q)^*$]

If  $K\in L^1([0,T], \Lin(\xR^n,\xR^n))$, i.e. $K=K(\cdot)$  is a $n\times n$-matrix-valued $L^1$ map on $[0,\bar S]$, we use the exponential notation  $[0,S]^2\ni (s_1,s_2) \mapsto e^{\int_{s_1}^{s_2} K}$ to denote the {\it fundamental matrix solution} of the (time-dependent) linear equation 
$\dot v = K v$: namely, for every $s_1,s_2\in [0,\bar S]$, $\bar v\in \xR^n$, $e^{\int_{s_1}^{s_2} K}\bar v = v(s_2)$, where $v(\cdot)$ is the solution to the linear Cauchy problem
$\dot v = K v,\,\,v(s_1) = \bar v$.

If $\mathcal{E}$ is a metric space with a metric $d$, and  $\bar e\in  \mathcal{E}$, ${ Q}\subseteq \mathcal{E}$, we use  $d ({\bar e},{ Q})$  to  denote {\it the distance of $\bar e$  from ${ Q}$}, namely $({\bar e},{Q}):=\inf\limits_{e\in{ Q}} d(\bar e,{e}).$

For any subset ${Q}\subseteq E$ of an Euclidean space $E$, by $\text{co}({Q})$ we mean the {\it convex hull of ${Q}$}, i.e. the smallest convex set containing ${Q}$, obtainable by intersection of all convex sets containing ${Q}$. The symbol $\overline{\text{co}}({Q})$ denotes the {\it closed  convex hull of ${Q}$}, i.e the closure of $\text{co}({Q})$, which happens to be smallest closed convex set containing ${Q}$.

For any subset ${Q}\subseteq X$ of a  topological space $X$, we use $\mathring{{Q}}$  to denote the set of interior points of ${Q}$.
If $n,m$ are positive integers,  $E\subseteq \rn$, is any subset and   $F:E \to \xR^m$ is a map,  differentiable  at some point $x\in \mathring{E}$, we use $DF(x)$ to denote the differential of $F$ at $x$.

If $H,K$ are sets, by using the notation $\mathcal{F}:H\setmap K$ we mean that $\mathcal{F}$ is a set-valued map from $H$ into $K$, that is a map from $H$ into the power set $\mathcal{P}(K)$.

To save space, with the expressions "Lipschitz function", "Lipschitz map", "Lipschitz vector field", we will mean 
"Lipschitz continuous  function" , "Lipschitz continuous  map", "Lipschitz continuous vector field", respectively.

{ If $A\subset \xR\times\xR^q$ is an open subset  and $F:A \to \xR^q $ is a time-dependent  vector field, continuous in $x$ and measurable in $t$, the expression {\it $x(\cdot)$ is a solution  of the differential equation $\dot x=F$ on an interval $I$} will always mean that $(t,x)\in A (t)$ for every $t\in I$ and  $x(\cdot)$ is  a {\it Carathéodory} solution of $\dot x=F(t,x)$, namely $x(\cdot)$ is absolutely continuous and  the equality  $\dot x(t)=F(t,x)(t)$ is verified at almost every $t\in I$.

More generally, if $X:A\setmap  \xR^q$ is a  set-valued  vector field ,  the expression {\it $x(\cdot)$ is  a solution of the  differential inclusion  $\dot x(t)\in X(t,x)$ on $I$} will  mean that $x(\cdot)$ is absolutely continuous   and  $\dot x(t) \in  X(t,x(t))$, for a.e. $t\in I$. }

\section{ Set-valued Lie  brackets  and  $QDQ$ Approximating Cones} 
\subsection{Lie brackets for Lipschitz vector fields}
If $f: \rn \to \rn$ is a Lipschitz vector field, we use $\Diff(f)$ to denote the set of   differentiability points of $f$, which,  by Rademacher's theorem, has full measure.    The following notion of  {\it set-valued Lie bracket  of two Lipschitz continuous vector fields $f,g: \rn \to \rn$} has been introduced in \cite{RampSuss2001}: $$[f,g]_{set}(x):=\overline{\text{co}}\Big\{v\in\rn,\,\,\,\,v=\lim\limits_{n\to \infty} [f,g](x_n),\,\,\,\, x_n \to x,\,\, (x_n)\subset \Diff(f)\cap \Diff(g)\Big\},$$
In this formula we mean that {\it all}  sequences 
$(x_n) \subset\Diff(f)\cap \Diff(g)$ are considered for which $x_n\to x$ and the limits $\lim\limits_{n\to \infty} [f,g](x_n)$ do exist.  Since $f,g,Df,Dg$ are bounded in a neighborhood of any point $x\in\rn$, by compactness one has $[f,g]_{set}(x)\neq \emptyset$. Moreover
$[f,g]_{set}(x)=\{[f,g](x)\}$ as soon as $f,g$ are of class $C^1$ in a neighborhood of $x$.
One trivially has that  the relations  $[f,f]_{set}=\{0\}$ and $[f,g]_{set}=-[g,f]_{set}$ keep holding for set-valued brackets, with the understanding that, for any subset $S$ of a vector space, $-S$ is the set of opposites of elements in $S$. Furthermore, some basic  results have been generalized to set-valued Lie brackets. For instance, the flow of $f,g$ locally commute if and only if $[f,g]_{set}=\{0\}$. Furthermore,  a Frobenius-type result  holds true for  Lipschitz distributions (see \cite{RampSuss2007}), {as well as} a local controllability theorem analogous to  Rashewski-Chow's \cite{RampSuss2001}.\footnote{{ In \cite{RampFeleqi}, a  notion of  iterated  set-valued Lie brackets  suitable for  local controllability issues  have been investigated as well.} }
\subsection{Quasi Differential Quotients  ($QDQ$s)}

{
	
	 We  call a  function   $\rho : \xR_+\to\xR_+\cup\{+\infty\}$   a   {\it pseudo-modulus}   if it  is monotonically nondecreasing and
	$\ds\lim_{s\to 0^+}\rho(s) =\rho(0)= 0$. 

Let us  recall the concept of  {\it Quasi Differential Quotient} for set-valued maps, which was  introduced  in \cite{RampPal} as a special case of Sussmann's Approximating Generalized Differential Quotient,  to address infimum gaps problems. 
\begin{dfntn}[Quasi Differential Quotients ($QDQ$)]\label{qdq}
	Let $N,n$ be non negative integers. Let $F : \xR^N \setmap \xR^n $ and  $\Lambda\subset \Lin\{ \xR^N, \xR^n\}$ be  a
	set-valued map and a compact subset, respectively, and let  $ \Gamma\subset\xR^N $  be any  subset. Given a pair $(\bar x,\bar y) \in \xR^N\times\xR^n  $,
	we
	say that $\Lambda$ is a {\em Quasi Differential Quotient  ($QDQ$) of $F$ at  $(\bar x,\bar y)$ in the direction of $\Gamma$} if there
	exists a pseudo-modulus $\rho$ such that for every  $\delta$ with $\rho(\delta)<+\infty$ there is a continuous map 
			$\displaystyle (L_\delta,h_\delta):\left(\bar{x}+B_{\delta}\right)\cap\Gamma\to \Lin( \xR^N, \xR^n) \times \xR^n$
			verifying
			\begin{equation}\label{InclusionedefQDQ}\begin{array}{c}\min_{L'\in\Lambda}|L_\delta(x) - L'|\leq \rho(\delta), 
					\quad |h_\delta(x)|\leq \delta \rho(\delta), \quad\hbox{and} \\[3mm]
					\bar y +  L_\delta(x)\cdot(x-\bar x)  + h_\delta(x)\in F(x),\end{array}\end{equation}
			whenever $x
			\in  (\bar x+B_\delta)\cap\Gamma$ .

\end{dfntn}\begin{rmrk} \label{QDQsemplificato}
	If $F$ is a single-valued, continuous map, one has necessarily $\bar y=F(\bar x)$, so the inclusion in \eqref{InclusionedefQDQ}  reduces to an equality. 
	%	Therefore, in order   to prove that $\Lambda$ is a $QDQ$ of $F$ at $(\bar \gamma,F(\bar \gamma))$ it is sufficient enough to find a modulus $\rho$ and  a family $\{L_\delta, \,\,\delta>0\}$ of continuous maps $L_\delta:\left(\bar{\gamma}+B_{\delta}\right)\cap\Gamma\to \Lin(\xR^N, \xR^n)$ satisfying $\ds\min_{L'\in\Lambda}|L_\delta(\gamma) - L'|\leq \rho(\delta), 
	%	\quad \text{ and} \quad|F(\gamma)-F(\bar \gamma)-L_\delta(\gamma)(\gamma-\bar{\gamma})|\leq \delta \rho(\delta) $
	%	whenever $\delta>0$ and $\,\gamma
	%	\in  (\bar\gamma+B_\delta)\cap\Gamma.$ 
	%%	Indeed, in this case the continuity of the error function \\ $\delta:=F(\gamma)-F(\bar \gamma)-L_\delta(\gamma)(\gamma-\bar{\gamma})$ would follow from the hypotheses.
	Furthermore,  if there exists  $\tilde{\varepsilon}>0$ and a continuous map $L:(\bar x+B_{\tilde{\varepsilon}})\cap \Gamma \to \Lin(\xR^N,\xR^n)$ satisfying $$\lim\limits_{\Gamma \ni x\to \bar x}\dist{L(x)}{\Lambda}=0,\quad \text{ and} \quad \lim\limits_{\Gamma \ni x\to \bar x} \frac{|F(x)-F(\bar x)-L(x)\cdot (x-\bar x)|}{|x-\bar x|}= 0,$$ then $\Lambda$ is a $QDQ$ of $F$ at $(\bar x,F(\bar x))$ in the direction of $\Gamma$. To see this, it is sufficient to define $L_\delta$ as the restriction of $L$ to $\bar x+ B_{\delta}\cap \Gamma$ and to  consider the  modulus $\rho:(0,\tilde{\varepsilon})\mapsto\xR$ defined by setting $$\rho(\delta):=\max\Bigg\{\sup\limits_{x \in B_\delta\cap \Gamma}\dist{L(x)}{\Lambda}, \sup\limits_{x \in B_\delta\cap \Gamma}\frac{|F(x)-F(0)-L(x)\cdot x|}{\delta}\Bigg\}\qquad \delta\in (0,\tilde{\varepsilon}). $$ 
\end{rmrk}

\subsubsection{Some basic properties of $QDQ$s}
\begin{prpstn}\cite{AngrisaniRampazzoQDQ} \label{basicprops} Let   
	$F,G:\xR^N \rightsquigarrow \xR^n$
	%	$F,G:\xR^N \rightsquigarrow \xR^n$%
	be set-valued maps. Assume that  $\bar x \in \xR^N$, $\bar y, \bar y_F, \bar y_G \in \xR^n$, $\Gamma, \Gamma_F, \Gamma_G \subseteq  \xR^N$,
	% $L \in \Lin(\xR^N,\xR^n)$,
	and $\alpha,\beta \in \xR$. Then: 
	\begin{enumerate}
		\item {\rm [Locality]} If $U$ is a neighborhood of $\bar x$ and $F(x)=G(x)$ for $x \in U\cap \Gamma$, then $\Lambda$ is a $QDQ$ of $F$ at $(\bar x, \bar y)$ in the direction of $\Gamma$ if and only if it is a $QDQ$ for $G$ at $(\bar x, \bar y)$ in the direction of $\Gamma$.
		\item {\rm [Linearity]} If $\Lambda_F$ and $\Lambda_G$ are $QDQ$ of $F$ and $G$ at points $(\bar x, \bar y_F)$ and $(\bar x, \bar y_G)$ in the direction of $\Gamma_F$ and $\Gamma_G$, respectively, then $\alpha\Lambda_F+\beta\Lambda_G$ is a $QDQ$ of $\alpha F+\beta G$ at point $(\bar x, \alpha\bar y_F+ \beta\bar y_G)$ in the direction of $\Gamma_F \cap \Gamma_G$.
		\item {\rm [Set product property]} Under the same assumptions as in $(2)$, $\Lambda_F\times \Lambda_G$ is a $QDQ$ at $\Big(\bar x, (\bar y_F, \bar y_G)\Big)$, in the direction of $\Gamma_F\cap \Gamma_G$, of the set-valued map $F \times G: x \rightsquigarrow F(x)\times G(x)$.
		\item {\rm [Product Rule]} If $m=1$, and still  using the same notation as in $(2)$, we have $F(\bar x)\Lambda_G+ G(\bar x)\Lambda_F$ is a $QDQ$ of $FG: x \rightsquigarrow F(x)G(x)$. \footnote{ We are using  $A+B$ to denote the set of sum of elements of two subsets $A,B\subset W$ of a vector space $W$. Furthermore, if $F$ is the field underlying a vector space $W$ and $\Omega\subset F$, we use the notation $\Omega A:=\{\omega a,\,\,\,\omega\in\Omega,  a\in A\}$. }
		\item \label{casospeciale} If $F$ is single-valued and $L\in  \Lin( \xR^N,\xR^n) $, $\{L\}$ is a $QDQ$ of $F$ at $(\bar x, \bar y)$ in the direction of $\xR^n$ if and only if $F$ is differentiable at $\bar x$ and $L=DF(\bar x)$.
	\end{enumerate}
\end{prpstn}
\begin{prpstn}[Chain rule]\label{chain}
	\cite{AngrisaniRampazzoQDQ}  	Let   $F:\xR^N \rightsquigarrow \xR^n$ and $G:\xR^n \rightsquigarrow \xR^l$ be set-valued maps, and consider the composition $G\circ F:\xR^N \ni x  \rightsquigarrow \bigcup\limits_{y \in F(x)} G(y) \in\xR^l$.
	Assume that $\Lambda_F$ is a $QDQ$ of $F$ at $(\bar x, \bar y)$ in the direction of $\Gamma_F$ and $\Lambda_G$ is a $QDQ$ of $G$ at  $(\bar y, \bar z)$ in a direction $\Gamma_G$ verifying $\Gamma_G\supseteq F(\Gamma_F)$. Then the set $\Lambda_G \circ \Lambda_F:=\Big\{ML,\quad M\in \Lambda_G,  \, L\in \Lambda_F\Big\}$
	%	 of all compositions of elements of $\Lambda_G$ with elements of $\Lambda_F$ 
	is a $QDQ$ of $G \circ F$ at $(\bar x, \bar z)$ in the direction of $\Gamma_F$.
\end{prpstn}
The following technical result will be essential for computing $QDQ$s of multiple variations of a given control process.
\begin{prpstn}\label{proppseudoaffine}
	Let $N,q$ be positive integers, and let  $F:\xR_+^N \to \xR^q$ be a map such that
	\begin{equation}\label{pseudoaffine}
		F(\bm{\varepsilon})-F(0)=\sum\limits_{i=1}^{N} \big(F(\varepsilon^i \mathbf{e}_i)- F(0)\big)+ o({|\bm\varepsilon|}),\qquad \qquad\quad\forall \bm{\varepsilon} =(\varepsilon^1,\dots,\varepsilon^N)\in\xR^N.\end{equation} 
	
	If,  for any $i=1,\ldots,N$,  $\Lambda_i\subset Lin(\xR,\xR^q)$ is a $QDQ$ at $0\in\xR$ of the map $\alpha\mapsto F(\alpha \mathbf{e}_i)$ in the direction of  a set $\Gamma_i\subseteq\xR_+$, then 
	the  compact set
	$$
	\Lambda :=\Big\{L\in Lin(\xR^N,\xR^q),\quad	L(v) = L_1v^1 + \ldots + L_N v^q,    \quad (L_1,\ldots, L_N)
	\in \Lambda_1\times\dots \times \Lambda_N\Big\}\qquad \footnote{As a set of matrices,  $\Lambda$ is made of $q \times N$ matrices whose $i$-th column is an element of $\Lambda_i$.}$$
	is a $QDQ$ of $F$ at $0$ in the direction of $\xR_+^N$.
	
\end{prpstn}

\begin{proof}
	Since, for every vector $k\in\xR^q$ and every map  $\psi:\xR_+^N\to \xR^q$ such that  $\psi (\bm{\varepsilon})=o(|\mathbf{\bm{\varepsilon}}|)$, the singleton $\{ 0\}\subset  \xR^N\times\xR^q$ is a $QDQ$ at $0$ of  the map $k+\psi$ in the direction of any subset of $\xR_+^N$,  in view of the linearity property  of the $QDQ$ it is not restrictive to assume the condition	$
	F(\bm{\varepsilon})=\sum\limits_{i=1}^{N} F(\varepsilon^i \mathbf{e}_i)$ instead of  \eqref{pseudoaffine}.\\ On the other hand, 
	since $ \varepsilon_i \mathbf{e}_i = P^i (\bm{\varepsilon})$ for every $i=1,\ldots,N$, where $P^i$ is the projection on the $i$-axis of $\xR^N$, 	$
	F(\bm{\varepsilon})=\sum\limits_{i=1}^{N} F(\varepsilon^i \mathbf{e}_i)$ reads
	$
	F=\sum\limits_{i=1}^{N} F\circ P^i.
	$
	Therefore the thesis follows from the chain rule, the linearity, and the  trivial fact that, for every $j=1,\dots,N$, the singleton $\{\mathbf{e}^j\}$ is a $QDQ$ (at any point and in any direction) of the projection $P^j$.
\end{proof}

We will regard the  {Clarke's Generalized Jacobian} as a particular $QDQ$. Let us begin recalling its definition:
\begin{dfntn}[Clarke's Generalized Jacobian]\label{CGJ}
	Let $F: \xR^N \to \xR^n $ be a map and assume that it is Lipschitz in a neighborhood  of a point $x\in \xR^N$. 
	The subset  $$\partial_x^C F(x):=\overline{\text{co}}\Big\{L=\lim\limits_{n\to \infty} DF(x_n),\,\,\, x_n \to x\,\,\,\, (x_n)\subset \Diff(F)\Big\}\subseteq \Lin(\xR^N,\xR^n)$$
	(where it is meant  that one takes limits along all sequences $(x_n)\subset \Diff(F) $  converging to $x$)	is called  {\it Clarke's Generalized Jacobian of $F$ at $x$.}
\end{dfntn}

\begin{prpstn}[Clarke's Generalized Jacobian is a $QDQ$]\label{Clarketh}
	Let  $\Omega\subset\xR^n$ be  an open set and let $f:\Omega\to\xR^n$  be a $K$-Lipschitz continuous map, for some $K>0$. 
	Then, for every $x^*\in\Omega$, the Clarke's Generalized Jacobian $\partial^C_x f(x^*)$ is a $QDQ$ of $f$ at $x^*$ in the direction of $\Omega$.
\end{prpstn}

We shall prove Proposition \ref{Clarketh} after stating the following technical result, whose proof is fairly straightforward.

 \begin{prpstn}\label{sucoleqdq3}
	Let  $x^* \in \xR^m$, and let $V$ and  $F:V\leadsto \xR^n$ be a compact neighborhood of $x^*$ and  a continuous set-valued map  with closed values, respectively. Let $\Gamma\subseteq\xR^m$ be a closed subset verifying $x^*\in\Gamma$, let  $\rho:[0,1]\to [0,+\infty[$ be a modulus,  and  let $y^* \in  F(x^*)$ and  $\Lambda\subseteq \Lin(\xR^m,\xR^n)$.
	Let $(x_j,y_j,F^j,\Lambda_j,\rho_j)_{j\in\xN}$ be a sequence such that:\begin{itemize}  
		\item[(i)] for every $j\in\xN$, 
		$\left (x_j,y_j\right)_{j\in \xN}\subset (\Gamma\cap V)\times\xR^n$, $F_i:V\leadsto \xR^n$ is a set-valued map, $\Lambda_j\in \Lin(\xR^m,\xR^n)$, and $\rho_j:[0,1]\to [0,+\infty[$ is a  monotonic non-increasing function;

		\item[(ii)]
		\bel{subs4}
		\begin{array}{c}\displaystyle
			\lim_{j\to\infty}\sup_{x\in V} d^\sharp\left(F_j(x),F(x)\right) =0\quad  \lim_{j\to\infty}
			d(y_j, y^*) = 0\\\\ \displaystyle
			\lim_{j\to\infty}  \sup_{\alpha\in[0,1]}|\rho_j(\alpha)-\rho(\alpha)| =  0 \quad \lim_{j \to\infty} d^\#(\Lambda_j,\Lambda) = 0;\end{array}\eeq

		\item[(iii)] for every $\delta\in]0,\bar\delta]$, where $\bar\delta>0$ is a suitable positive number verifying $B(x^*,\bar\delta)\subset V$,  there exists a sequence \begin{center} $\left(\Big({L_\delta^j},{h_\delta^j}\Big)\right)_{j\in\xN}\subset C^0\Big(\bar{x}+B_{\delta},\Lin\{ \xR^m, \xR^m\} \times \xR^m  \Big)$  and $ \Big({L_\delta},{h_\delta}\Big)\in C^0\Big(\bar{x}+B_{\delta},\Lin\{ \xR^m, \xR^m\} \times \xR^m  \Big)$ \end{center}
		such that
		$$
		\Big({L_\delta^j},{h_\delta^j}\Big)\to \Big({L_\delta},{h_\delta}\Big)
		$$
		uniformly, and, for every  $j\in\xN$ and  $x
		\in  (x_j+B_\delta)\cap\Gamma$, verifies
		\bel{qqdq}\begin{array}{c} 
			y_j +  {L_\delta^j}(x)\cdot(x-x_j)  + {h_\delta^j}(x)\in F^j(x)\\\\ 
			\min_{L'\in\Lambda_i}|{L_\delta^j}(x) - L'|\leq \rho_j(\delta), 
			\quad |{h_\delta^j}(x)|\leq \delta \rho_j(\delta);\end{array}\eeq

	\end{itemize}
	Then  $\Lambda$ is a $QDQ$ of $F$ at $(x^*,y^*)$ in the direction of  $\Gamma$. More precisely, 
	for every $\delta\in [0,\frac{\bar\delta}{2}]$,  the map
	$(L_\delta,h_\delta)$ 
	verifies
	\bel{qdq}\min_{L'\in\Lambda}|L_\delta(x) - L'|\leq 2\rho(\delta), 
	\quad |h_\delta(x)|\leq 2\delta \rho(\delta), \quad\hbox{and}\
	y^* +  L_\delta(x)\cdot(x-x^*)  + h_\delta(x)\in F(x),\eeq
	whenever $x
	\in  (x^*+B_\delta)\cap\Gamma$ .
\end{prpstn}

\begin{proof}[Proof of Proposition \ref{Clarketh}]
	Let us consider the standard mollifier $\displaystyle
	\eta(x) := \frac{1}{I_n}\exp\left(-\frac{1}{(1-|x|^2)}\right) {\bf 1}_{B(0,1)} $
	where $I_n$ is the normalizing factor  $I_n:=\displaystyle\int_{B(0,1)} \exp\left(-\frac{1}{(1-|x|^2)}\right) dx$,
	and, for every $\epsilon>0$, let us set
	$\displaystyle \eta_\epsilon(x):=\frac{1}{\epsilon^n} \eta\left(\frac{x}{\epsilon}\right)$.
	For every $x\in \Omega$ and $\epsilon>0$ so small that $B(0,\epsilon)\in \Omega$, let us {consider the mollification}
	$$
	f^\epsilon(x):= \eta_\epsilon*f (x)= \int_{B(0,\epsilon)} f(x-y)\eta_\epsilon(y){ dy }  =\int_\Omega f(x-y)\eta_\epsilon(y){ dy } .
	$$
	In particular,  for every $\epsilon>0$ one has
	$
	\ds\sup_{x\in \Omega}|f^\epsilon(x)-f(x)|\leq K\epsilon.
	$
	Moreover, for every $x\in \Omega$  and every $\epsilon<d(x,\partial\Omega)$, one has
	\bel{derivata}
	Df^\epsilon(x)= \eta_\epsilon*Df (x)\in co\Big\{Df(y),\,\,\,y\in \text{Diff}(f), y\in B(x,\epsilon)\Big\}.
	\eeq
	Let us consider  the sequence of (single-valued) maps $(F^j)_{j\in\xN}$ defined by setting $\ds F^j:=\displaystyle f^{\frac{1}{j^2}}$ for every $j\in\xN$. Since the maps  $F^j$ are  smooth,
	for every $j\in\xN$ such that $\ds d(x^*,\partial\Omega) > \frac 1j+\frac{1}{j^2}$ and every 
	$x\in B\left(x^*,\ds \frac 1j\right)$, one has 
	\bel{ap1}
	F^j(x) = F^j(x^*) +L^{j}(x)\cdot (x-x^*),\eeq where we have set $$\displaystyle L^{i}(x): = 
	\int_0^1 D F^j(x^*+t(x-x^*))dt = \int_0^1 \left( \int_\Omega \eta_{\frac{1}{j^2}}(y)D F(x^*+t(x-x^*)-y)dy\right)dt .$$
	Aiming to apply Proposition \ref{sucoleqdq3}, let us  set, for every $j\in\xN$ such that $\ds d(x^*,\partial\Omega) > \frac 1j+\frac{1}{j^2}$  and every $\delta>0$, $$\begin{array}{c}\Lambda^j := co \left\{Df(y),\,\,\,y\in Diff(f), y\in B\left(x^*, \ds \frac 1j+\frac{1}{j^2}\right)\right\}\qquad \Lambda:=\partial_x^Cf(x^*), \\\\  y_j:=F^j(x^*),\quad L_\delta^j:= L^{\left[\frac {1}{\delta}\right]}\\\\ h^j_\delta(x):= \left(F^j(x)-F^{\left[\frac {1}{\delta}\right]}(x)\right) -\left(F^j(x^*)-F^{\left[\frac {1}{\delta}\right]}(x^*)\right) = \left(L^j(x)-L^{\left[\frac {1}{\delta}\right]}(x)\right)\cdot (x-x^*) .
	\end{array}$$
	Observe that, for every $j\in\xN$ and $\ds x\in  B\left(x^*, \frac 1j\right)$, we have 
	$\ds
	L^j(x)\in \Lambda^j.
	$
	Setting 
	$$
	\rho_j(\delta) := \frac1\delta \sup_{x\in B(x^*,\delta)}  \left|\left(L^j(x)-L^{\left[\frac {1}{\delta}\right]}(x)\right)\cdot (x-x^*)\right|
	$$ 
	one gets (that  $\rho_j$ is  monotonically decreasing and) $\rho_j(\delta) \leq 2K\left(\delta + \frac{1}{\delta j^2}\right)$ for every $j\in \xN$ and $\delta>0$, so that the pointwise limit $\ds\rho(\delta)=\lim_{k\to\infty}\rho_j(\delta)$ is a modulus, i.e.  $\rho(\delta)$ (is monotonically deceasing and)  verifies $\ds\lim_{\delta\to 0}\rho(\delta) = 0$.

	Since, for every $x\in B(x^*,\delta)$,
	$$\begin{array}{c}
		F^j(x) = y_j + L^{j}_\delta (x)\cdot (x-x^*) +  h^j_\delta(x), \\\\
		\ds\min_{L'\in\Lambda_i}| L^{j}_\delta (x)-L'| \leq  | L^{j}_\delta (x)-L^j|\leq     \rho_j(\delta)\qquad    h^j_\delta(x)
		\leq 2\left(\frac{1}{j^2} + \delta^2\right) = \rho_j(\delta)\delta,
	\end{array} 
	$$
	in view of  Proposition \ref{sucoleqdq3} the proof is concluded provided  
	$\lim_{j \to\infty} d^\#(\Lambda_j,\Lambda) = 0.$ \\ Actually, this is easily verified: indeed, if
	there existed $a>0$ such that there were $M^j\in \Lambda^j $ such that
	$d(M^j,\partial^C_xf(x^*))>a$  for all $j\in\xN$ , there would be a subsequence $(M^{j_k})_{k\in\xN}$ of $(M^j)_{j\in\xN})$ converging to  an element $M\in \Lin(\xR^m,\xR^n)$.\footnote{Let us point out that for all $j\in\xN$ and all $L\in \Lambda^j$, one has  $|L|\leq K$, where $K$ is the Lipschitz constant of $f$.} This would imply that $M\in \partial_x^Cf(x^*)$, so contradicting the inequality $d(M,\partial_x^C f(x^*))\geq a(>0)$. \end{proof}

The following results says that, for every vector field measurable in time and Lipschitz in space,  {the set-valued image at $t$ of the corresponding variational differential inclusion is a $QDQ$ of the flow map, for  any $t$ in the interval of existence.} 
%the set-valued map of the solutions of the variational differential inclusion is a $QDQ$ of the flow map.

\begin{lmm}\label{variazionale}

	Let $S>0$ and $F:\xR\times\xR^n\to \xR^n$ be a vector field such that, for every $s\in\xR$, $F(s,\cdot)$ is Lipschitz with Lipschitz constant $L_F(s)$ satisfying $\int_0^S L_F(s)\,ds=A<\infty$, and, for every $y\in\xR^n$, $F(\cdot,y)$ is a bounded Lebesgue measurable map. Let $q\in\xR^n $, and let us assume that the Cauchy problem $$\begin{cases} y'(s)=F(s,y(s))\\ y(0)=\xi \end{cases}$$
	has a {(necessarily unique)} solution $s\mapsto \Phi_s^F(\xi)$ on an interval $[0,S]$ for every $\xi$ in a neighborhood $U$ of $q$.
	Then, for every  $s\in [0,S]$
	the set,
	$$\Lambda=\left\{L(s),\,\, \text{$L$ is a solution of}\,\,\,L'(\sigma) \in \partial_y^C F(\sigma, y(\sigma))\cdot L(\sigma),\,\, L(0)=\mathbf{1}\right\}$$
	is a $QDQ$ of the map  $\xi\mapsto \Phi_s^F(\xi)$ at $q$ in the direction of $\xR^n$.
\end{lmm}
{ \begin{proof}
		Let $\eta:\xR^n \to \xR_+$ be  the standard mollifier as in Proposition \ref{Clarketh}, and for any $\sigma>0$, let us set  $\eta_\sigma(y)=\frac{1}{\sigma^n}\eta\left(\frac{y}{\sigma}\right)$.
		For every $t\in\xR$, let us consider the convolution $\ds \xR^n\ni y\mapsto F_\sigma(t,y):=\int_{\xR^n} F(t,y-h)\eta_\sigma(h) dh$.
		Since  the vector field $F_\sigma$ is   $C^\infty$ with respect to $y$ and measurable with respect to $t$, there exists a unique solution  $y_\sigma(s)$ to the Cauchy problem on $[0,S]$ $$\begin{cases} y'(s)=F_\sigma(s,y(s))\\ y(0)=\xi. \end{cases}$$
		It is easy to  check that $y_\sigma(s)$ uniformly converges to $y(s)$.	Indeed, \begin{multline*}|y_\sigma(s)-y(s)|\le \int_{0}^s \left|F_\sigma(\tau,y_\sigma(\tau))-F(\tau,y(\tau))\right|\,d\tau \le \\ \le \int_0^s \left| F_\sigma(\tau,y_\sigma(\tau))-F(\tau,y_\sigma(\tau)) \right|\,d\tau + \int_0^s \left| F(\tau,y_\sigma(\tau))-F(\tau,y(\tau)) \right|\,d\tau \le \\ \le \int_0^s \int_{\xR^n}\left| F(\tau,y_\sigma(\tau)-h)-F(\tau,y_\sigma(\tau))\right|\eta_\sigma(h)\,dh\,d\tau+\int_0^s L_F(\tau)|y_\sigma(\tau)-y(\tau)| \,d\tau \le \\ \le  \int_0^s L_F(\tau)\,d\tau \cdot \int_{\xR^n}h\eta_\sigma(h)\,dh+\int_0^s L_F(\tau)|y_\sigma(\tau)-y(\tau)|\,d\tau\le A\sigma+\int_0^s L_F(\tau)|y_\sigma(\tau)-y(\tau)|\,d\tau, \end{multline*}
		so that,	using Gronwall's Lemma, we get $$|y_\sigma(s)-y(s)| \le  A\sigma+A\sigma\int_0^S L_F(\tau)e^{\int_0^S L_F(r)\,dr}\,d\tau \le (A+A^2e^A)\sigma.$$  In addition, by classic theory of ODE's we have, for every $\varepsilon\in\xR^n$, \begin{equation}\label{classicvariational}\Phi_s^{F_\sigma}(q+\varepsilon)=\Phi_s^{F_\sigma}(q)+M_\sigma(s)\cdot \varepsilon+ o(|\varepsilon|),\end{equation}
		where the matrix-valued map $M_\sigma$ is the solution of the variational Cauchy problem $$\begin{cases}\ds M'(s)=\frac{\partial F_\sigma}{\partial y}(s,y_\sigma(s))\cdot M(s),\\M(0)=\mathbf{1}. \end{cases}$$ Now, since $F(s,\cdot)$ is Lipschitz, the spatial gradient of the mollified function $F_\sigma$  coincides with the mollification of  the ($L^\infty$) gradient $\ds \frac{\partial F}{\partial y}(s,\cdot)$, i.e.  $\ds\frac{\partial F_\sigma}{\partial y}(s,y)=\int_{\xR^n} \frac{\partial F}{\partial y}(s,y-h)\eta_\sigma(h)\,dh.$ Therefore, with computations similar to those performed above, it follows that
		
		$$d(M_\sigma(S),\Lambda)\to 0 \qquad \text{ as } \qquad \sigma \to 0.$$ Finally, if {$\displaystyle\lim_{\varepsilon \to 0}\sigma(|\varepsilon|)=0$}, we obtain \begin{multline}\Phi_S^F(q+\varepsilon)=\Phi_S^{F_\sigma}(q+\varepsilon)+o(|\varepsilon|)=\Phi_S^{F_\sigma}(q)+M_\sigma(S)\cdot \varepsilon+ o(|\varepsilon|)=\Phi_S^F(q)+M_\sigma(S)\cdot \varepsilon+o(|\varepsilon|)\end{multline} with $M_\sigma(S)$ having vanishing distance from $\Lambda$ as $|\varepsilon|\to 0$. From this, our thesis follows according to Remark \ref{QDQsemplificato}.
\end{proof}}\subsection{$QDQ$-approximating cones and multicones}
Let  $V$ be  a finite-dimensional real vector space. A subset $C\subseteq V$ is called a {\it cone} if $\alpha v\in C, \forall \alpha\ge 0$ and $\forall v \in C.$ A family $\mathcal{C}$  whose elements are cones is called a {\it multicone}. A {\it convex multicone} is  a multicone whose elements are convex cones.

For any given subset $E \subseteq V$, the set $E^\perp:=\{v \in \xR^n, \, v \cdot c \le 0 \,\,\forall c \in C\}\subseteq E^*$ is a closed cone, called the {\it polar cone of $E$}.

Let us introduce the notion of transversality, according to \cite{Sussmann2}.
Two cones $C_1$ and $C_2$ are said to be {\it transversal} if $C_1-C_2=V$, where we use the notation $\ds C_1-C_2:=\big\{c_1-c_2, \quad (c_1,c_2)\in C_1\times C_2\big\}$. Two multicones $\mathcal{C}_1$ and $\mathcal{C}_2$ called {\it transversal} if  $C_1\in \mathcal{C}_1$ and $C_2\in \mathcal{C}_2$ are transversal as soon as$(C_1,C_2)\in \mathcal{C}_1\times\mathcal{C}_2$.

Two transversal cones $C_1$ and $C_2$ are called {\it strongly transversal} if $C_1\cap C_2 \supsetneq \{0\} $. This is trivially equivalent to the existence of a non-zero linear form $\mu$  and an element $c\in C_1\cap C_2$  such that $\mu c >0$. More generally, we say that two transversal multicones $\mathcal{C}_1$, $\mathcal{C}_2$ are {\it strongly transversal} if there exists a non-zero linear form $\mu$ such that for any choice of cones $C_i \in \mathcal{C}_i$, $i=1,2$, there is an element $c\in C_1\cap C_2$ such that $\mu c >0$.
One says that two cones $C_1$, $C_2$ are {\it linearly separated if  $C_1^\perp \cap -C_2^{\perp}\supsetneq \{0\}$}, namely there exists a form $\mu\in V^*\backslash\{0\}$ such that $\mu c_1\geq 0, \mu c_2\leq 0$ for all $(c_1,c_2)\in C_1\times C_2$. It is easy to check that $C_1$ and $C_2$ are linearly separated if and only if they are not transversal. {For multicones one has the following fact:}
\begin{lmm} \label{nontrasversalitaforte} \cite{Sussmann2}  
	Let $\mathcal{C}_1$ and $\mathcal{C}_2$ be two multicones that are not strongly transversal. If there is a linear functional $\mu$ that is in $C_2^\perp$ but not in $-C_2^\perp$ for all $C_2 \in \mathcal{C}_2$ , then there are two cones $C_1 \in \mathcal{C}_1$ and $C_2 \in \mathcal{C}_2$ that are not transversal, i.e.  $C_1, C_2$ are linearly separated.
\end{lmm}
{ \begin{dfntn}[$QDQ$-approximating multicones]\label{qdqmulticones}
		Let $E$ be any subset of an Euclidean space $\rn$ and $x \in E$\footnote{The definition of approximating multicone can be easily generalized to the case of  subsets of a differential manifold \cite{RampPal}.}. A convex multicone $\mathcal{C}$ is said to be a {\it $QDQ$-approximating multicone to $E$ at $x$} if there exists a  set-valued map $F:\xR^N \setmap \xR^n$, a convex cone $\Gamma\subset\xR^N$, and a $QDQ$  $\Lambda$ of $F$ at $(0,x)$ in the direction of $\Gamma$   such that  $$ F(\Gamma)\subseteq E, \qquad\mathcal{C}=\{L\cdot \Gamma, \,\, L \in \Lambda\}.$$
		When $\Lambda$ is a singleton, i.e. $\Lambda=\{L\}$, one  simply says  that {\it $C:= L\cdot \Gamma$ is a $QDQ$-approximating cone to $E$ at $x$.}
\end{dfntn}}
\begin{dfntn}[Local separation of sets]
	Two subsets $E_1$ and $E_2$ are {\it locally separated at $x$} if  there exists a neighborhood $U$ of $x$ such that $\displaystyle E_1\cap E_2 \cap U =\{x\}.$
\end{dfntn}
{ As a consequence of an open mapping theorem, the following fact holds true (see Theorem 4.37, p. 265 in \cite{Sussmann} where the lemma was proven in the more general context of $AGDQ$'s, of which $QDQ$ are a special case)

	\begin{lmm}\label{OpenMappingConsequence}{If two subsets $E_1$ and $E_2$ are locally separated at $x$ and if $\mathcal{C}_1$ and $\mathcal{C}_2$ are  $QDQ$-approximating multicones  for $E_1$ and $E_2$, respectively, at $x$, then  $\mathcal{C}_1$ and $\mathcal{C}_2$ are not strongly transverse.}
\end{lmm}}

\section{The minimum problem and the main result}

\subsection{The minimum problem}
The optimal control problem we are going to address, which will still label  $(P),$ is more general  than the one presented in the Introduction, {in that it involves a Lagrangiean $l$ as well as  an additional, bounded control $a$}.  Precisely we will consider the problem 
$$(P)\qquad
\left.
\begin{array}{l}
	\quad \ds	
	\min_{u\in\mathcal{U}} \left(\Psi(T,x(T))+\displaystyle \int_0^T l(x(t),u(t),a(t))\,dt\right),\\ \\
	\begin{cases}
		
		\displaystyle \frac{dx}{dt}=f(x,a)+\sum\limits_{i=1}^m g_i(x)u^i, \quad \text{ a.e. } t \in [0,T],\\
		\displaystyle \frac{d\nu}{dt}=|u|,\\
		\displaystyle (x,\nu) =(\hat{x},0),\qquad
	\end{cases} \qquad (T,x(T),\nu(T)) \in \mathfrak{T}\times [0,K]\end{array}\right.$$

where:\begin{itemize}\setlength\itemsep{-0.4em} \item[i)] the state variable $x$ belongs to $\xR^n$, for some $n>0$; \item[ii)] the vector fields  $g_i:\xR^n \to \xR^n, \, i=1,\ldots,m$  are locally Lipschitz; \item[iii)] the unbounded controls $u=(u^1,\ldots,u^m)$   take values in a  closed cone $C=C_1\times C_2$, where, for some non negative integers  $m_1$ and $m_2$  such that $m=m_1+m_2$,  $C_1 \subseteq \xR^{m_1}$ is a closed cone containing the coordinate axes, and $C_2 \subseteq\xR^{m_2}$ is a closed cone which does not contain any  straight line; the control $a$ takes values in a compact set $A\subset\xR^q$; \item[v)] the {\it drift} $f$ is continuous  and, for any value of the control $a\in A$, the function $x \mapsto f(x,a):\xR^n \to \xR^n$  is locally Lipschitz;
	\item[iv)] the real-valued {\it Lagrangian} $l:=l(x,u,a)$ has the form $l(x,u,a)=l_0(x,a)+l_1(x,u)$ and  is  continuous; furthermore, the map $x\mapsto l(x,u,a)$ is locally Lipschitz, uniformly for every $(u,a)\in C\times A$; moreover,   the so-called {\it recession function}  \begin{equation}\label{recfunc}
\ds\hat{l}_1(x,w^0,w):=\lim\limits_{r\to w^0} rl_1\left(x,\frac{w}{r}\right)\end{equation} is well-defined and locally Lipschitz with respect to $x$, uniformly as $(w^0,w)$ ranges on the bounded set  $[0,1]\times(C\cap B_1)$ ;
	\medskip
	\item[v)]  the {\it final cost} $\Psi(t,x)$ is Lipschitz, $0 \le K\le +\infty$,  the (time-dependent) {\it target} $\mathfrak{T}\subseteq \xR_{+}\times \xR^n$ is  a closed subset, and the  $L^1$ bound $\|u\|_1\leq K$ has been written in the equivalent  form $\nu\le K$.; \item[iv)] the minimization is performed over the set of the {\it  strict sense feasible processes}, where by {\it  strict sense process} \footnote{We use the expressions 'strict sense' in order to distinguish processes and controls  of the original  problem from those of the extended problem we introduce later, which  will be named 'extended sense processes'  and 'extended sense controls', respectively} we mean a six-tuple  $(T,u,a,x,\nu)$
	such that $(T,u,a) $  belongs to the family of  {\it strict sense controls}  $$\mathcal{U}:=\bigcup\limits_{T>0} \{T\} \times L^1([0,T],C\times A)$$ and 
	$(x,\nu)$ is the solution of the above control  system, whereas a strict sense   process is called  {\it feasible}  as soon as $(T,x(T),\nu(T)) \in \mathfrak{T}\times [0,K]$.

\end{itemize}

\begin{rmrk}\label{manifold}
	All the involved objects having an intrinsic character, the optimal control  problem and the corresponding results presented in this paper  can be easily extended   to the more general situation where the state $x$ range over an $n$-dimensional manifold.
\end{rmrk}

\begin{rmrk}\label{Timedependence}
	The  assumption that the Lagrangian cost $l$, the drift term $f$, and the vector fields $g_i$ are not time-varying can be removed by the standard procedure of regarding the time variable $t$ as an extra-state variable $x^0$ subject to the trivial equation $\dot x^0=1$
\end{rmrk}

\begin{dfntn} We say that $(\overline{T},\overline{u},\overline{a},\overline{x},\overline{\nu})$ is a strict sense {\it weak local minimizer} for problem (P) if there exists $\delta>0$ such that $$\Psi(\overline{T},\overline{x}(\overline{T}))+\int_0^{\ol{T}}l(\ol{x}(t),\ol{u}(t),\ol{a}(t))\,dt \le \Psi(T,x(T))+\int_0^{T}l(x(t),u(t),a(t))\,dt$$
	for all feasible processes $(T,u,a,x,\nu)$ such that $|T-\overline{T}|+\|(x,\nu)-(\ol{x},\ol{\nu})\|_{\infty} + \|(u,a) -(\bar u,\bar a) \|_1
	<\delta.$  {Actually, since $(x,\nu,u,a)$ and $(\ol{x},\ol{\nu},\ol{u},\ol{a})$ may have different domains, we tacitly  extend   $(x,\nu)$ and $(\ol{x},\ol{\nu})$ continuously from $[0,T]$ and $[0,\ol{T}]$ to $\xR_+$ so  that they are constant on $[T,+\infty]$ and $[\ol{T},+\infty]$. Furthermore, we extend  $(u,a)$ and $(\ol{u},\ol{a})$ by setting  $(u,a)(t)=(0,\hat{a})$ for any $t>T$ and $(\ol{u},\ol{a})(t)=(0,\hat{a})$ for any $t>\ol{T}$, for  some (irrelevant) choice of  $\hat{a}\in A$.}
	
\end{dfntn}

\subsection{The extended problem}
Since we are interested in necessary conditions  for minima, in principle we might ignore the existence problem. Yet, the unboundedness of the controls and the lack of adequate coercivity assumptions make the existence of an optimal control a quite unlikely{,   if not impossible. To get existence of {\it minima} it is then convenient  to continuously  embed   the problem in a more general one where  trajectories are somehow allowed to evolve also in a  degenerate time interval (consisting of a single time instant).  A distributional embedding being ruled out  because of the non-linearity of the problem,  a robust extension consists instead in first transforming   the original problem in a time-space problem, where the trajectories are replaced by their graphs, and, secondly,  considering the $C^0$-closure of the set of such graphs as the new minimization domain.} Precisely we will consider the extended problem 
$$(P_{ext})\quad \left.\begin{array}{l}
	\quad \ds	 
	\min_{(S,w^0,w,\alpha)\in\mathcal{W}} \left(\Psi(y^0(S),y(S))+\displaystyle\int_0^S l^e((y,w^0,w,\alpha)(s))\,ds\right),\\ \\
	\begin{cases} \displaystyle \frac{dy^0}{ds}(s)=w^0(s),\\ \displaystyle \frac{dy}{ds}(s)=f(y(s),\alpha(s))w^0(s)+\sum\limits_{i=1}^m g_i(y(s))w^i(s),\\ \displaystyle \frac{d\beta}{ds}(s)=|w(s)|,\\ \displaystyle (y^0,y,\beta)(0)=(0,\hat{x},0). \end{cases}\qquad (y^0(S),y(S),\beta(S)) \in \mathfrak{T}\times [0,K]\end{array}\right.$$

where:\begin{itemize} 
	\item[i)] $t=y^0\in\xR_+$ stands for the time parameter and $y\in\xR^n$ denotes the state variable;
	\item[ii)]  the {\it extended Lagrangian} $l^e$ is  defined by setting  $$l^e(x,w^0,w,\alpha):=l_0(x,\alpha)w^0+\hat{l}_1(x,w^0,w)\quad\forall (x,w^0,w,\alpha)\in \xR^n \times  \xR_+\times C\times A, $$
	$\hat l_1$ {being the recession function defined in \eqref{recfunc}.}\footnote{ In view of the sublinearity of $l$ in $u$, 
		$l^e$ is well-defined. As an example, one can consider the Lagrangian $l(x,u,a)=l_0(x,a)+\ell(x) |u|^r$ for some $r\in[0,1]$ and some Lipschitz function $\ell$, in which case  one  has $l^e(x,w^0,w,\alpha)=l_0(x,\alpha)w^0+\ell(x) |w|^r (w^0)^{1-r}$.}  \item[iii)]  the four-tuples  $(S,w^0,w,\alpha)$ belong to the set  $$\mathcal{W}:=\bigcup\limits_{S>0} \{S\} \times \Big\{(w^0,w,\alpha) \in L^\infty([0,S],\xR_+\times C\times A): \essinf(w^0+|w|)>0\Big\},$$ whose elements are called {\it extended sense controls}.; \item[iv)] the minimization is performed over the set of  {\it  extended sense  feasible processes}, which are defined as follows: \begin{itemize} \item  an {\it  extended sense   process}   is  a   seven-tuple  $(S,w^0,w,\alpha,y^0,y,\beta)$
		such that $(S,w^0,w,\alpha) $   is an  {\it extended sense control} and $(y^0,y,\beta)$ is the corresponding solution of the extended  Cauchy problem  in $(P)_{ext}$, whereas
		\item  an extended sense    process $(S,w^0,w,\alpha,y^0,y,\beta)$ is called  {\it feasible}  as soon as $(y^0(S),y(S),\beta(S)) \in \xR_*\times \mathfrak{T}\times [0,K]$.
	\end{itemize}

\end{itemize}

Let us give the notion local minimizer  for the  extended problem:
\begin{dfntn}  We say that $(\overline{S},\overline{w}^0,\overline{w}, \ol{\alpha},\ol{y}^0,\ol{y},\overline{\beta})$ is a {\it weak local minimizer} for problem $(P)_{ext}$ if there exists $\delta>0$ such that \begin{multline*}\Psi(\ol{y}^0(\overline{S}),\ol{y}(\overline{S}))+\int_0^{\ol{S}}l^e(\ol{y}(s),\overline{w}^0(s),\overline{w}(s),\ol{\alpha}(s))\,ds \le\\ \le \Psi({y}^0({S}),{y}({S}))+\int_0^{{S}}l^e({y}(s),{w}^0(s),{w}(s),{\alpha}(s))\,ds\end{multline*}
	for all feasible processes  $(S,w^0,w,\alpha,y^0,y,\beta)$ such that $$|S-\overline{S}|+\|(y^0,y,\beta)-(\ol{y}^0,\ol{y},\overline{\beta})\|_{\infty}+\|(w^0,w,\alpha)-(\ol{w}^0,\ol{w},\ol{\alpha}))\|_{1}
	<\delta.\quad  \footnote{ Since $(y^0,y,\beta,w^0,w,\alpha)$ and $(\ol{y}^0,\ol{y},\ol{\beta},\ol{w}^0,\ol{w},\ol{\alpha})$ may have different domains, as before we tacitly extend them to $\xR_+$ in the following way: the trajectories $(y^0,y,\beta)$ and $(\ol{y}^0,\ol{y},\ol{\beta})$  are prolonged to the whole $[0,+\infty[$ in such a way they are  continuous  and constant on $[S,+\infty[ $ and  $[\ol{S},+\infty[$, respectively; as for the controls,  the extensions consist in choosing a value $\hat{a}\in A$ and assigning the common constant value  $(0,\bf 0,\hat{a})$  to  both  $(w^0,w,\alpha)$ and $(\ol{w}^0,\ol{w},{\hat\alpha})$ on the intervals $]S,+\infty[ $ and  $]\ol{S},+\infty[$, respectively.}$$
\end{dfntn}
There is an obvious one-to-one {correspondence} between  strict sense processes and extended sense processes such that $w_0(s)>0$ for almost any $s\in (0,S)$. More precisely,
{\it A six-tuple $(T,u,a,x,\nu)$ is a strict sense process if and only if for every $S>0$ and every  strictly increasing, Lipschitz continuous, surjective function $\tau:[0,S]\to [0,T]$, the seven-tuple 
	\begin{multline*}(S,w^0,w,\alpha,y^0,y,\beta)(s):=\Big(S,w^0(s), u(\tau(s)) w^0(s) ,a(\tau(s)),\tau(s),x(\tau(s)),\nu(\tau(s))\Big)\qquad \forall s\in [0,S],\\ w^0(s):=\frac{d\tau}{ds}(s),\quad \text{for a.e.} s\in [0,S] \end{multline*}
	is an extended sense process with $w^0>0$ a.e. and, moreover, 
	
	$$\Psi(T,x(T))+\int_0^{T}l(x(t),u(t),a(t))\,dt = \Psi({y}^0({S}),{y}({S}))+\int_0^{{S}}l^e({y}(s),{w}^0(s),{w}(s),{\alpha}(s))\,ds.  $$
	
} This {one-to-one} correspondence preserves feasibility of a process,  and minima of the strict sense problem correspond  to minima for the restriction of the   space-time extended problem to processes having $w^0(s)$ almost everywhere positive (see \cite{ArMotRamp}).

\begin{rmrk}
	An important feature of the extended system is its  \textit{rate-independence}. By {this} we mean that if $ \sigma:[0,\hat S]\to [0,S]$ is a bi-Lipschitz function , then $(S,w^0,w,\alpha,y^0,y,\beta)$ is a feasible process if and only if $\Big(\hat S,(w^0\circ\sigma) \cdot \frac{d\sigma}{ds} ,(w\circ\sigma) \cdot \frac{d\sigma}{ds},\alpha\circ\sigma,y^0\circ\sigma,y\circ\sigma,\beta\circ\sigma\Big)$ is a feasible process. Two processes obtained one from the other in this way are called \textit{equivalent} and it is straightforward to verify they have the exact same associated costs, so that being a weak local minimizer is a property shared by equivalent processes. As it was observed in \cite{ArMotRamp}, this rate-independence implies the following fact:
	
{ 	\begin{itemize}
		\item  It is not restrictive to assume a minimizer  to be {\it canonical}, meaning that  $w^0(s)+|w(s)|= 1$ for almost every $s\in [0,\overline S].$
	\end{itemize}}
\end{rmrk}

To save space, for all $(y,w^0,w,a)\in \xR^{n}\times \xR_+\times C\times A  $ let us  introduce  the notation $$F^e(y,w^0,w,a):= f({y},a){w}^0+\sum\limits_{i=1}^m g_i({y}){w}^i$$ and 
$$
\cl{F}(y,w^0,w,a):=\left(w^0,F^e(y,w^0,w,a),l^e(y,w^0,w,a)\right).$$

Let us  define the Hamiltonian $H$ by setting, for every $(y,p_0,p,\lambda,\pi, w^0,w,a) \in $ $\xR^n\times (\xR^{1+n+1+1})^*\times W\times A $, $$H(y,p_0,p,\lambda,\pi, w^0,w,a):=p_0w^0+pF^e(y,w^0,w,a)-\lambda l^e(y,w^0,w,a)+\pi|w|.$$

We are now in the position of stating our main result:

\begin{thrm}[\bf A "higher order"  Maximum Principle]\label{TeoremaPrincipale}
	Let $\opprocess$ be a canonical local minimizer for the extended  problem  $(P_{ext})$, and
	%, with our hypotheses on $g_i,f,l,\hat{l}_1,\Psi,\mathfrak{T},A$ and $C$. 
	let $\mathbcal{T}$ be a $QDQ$-approximating multicone to the target set $\mathfrak{T}$ at $(\ol{y}^0,\ol{y})(\ol{S})$.

	Then there exist multipliers $(p_0,p,\lambda,\pi) \in \xR^*\times AC\Big([0,\ol{S}];(\xR^n)^*\Big) \times \mathbb{R^*} \times \mathbb{R^*}$ %and $\pi\le 0$
	such that $\pi\le 0$ (with  $\pi=0$ as soon as $\|\ol{w}\|_1<K$), $\lambda\geq 0$ and the following conditions are satisfied:\begin{itemize}
		\item[\rm\bf i)]{\sc(non triviality)}\,\,\, $(p_0,p,\lambda)\neq 0;$ \,\,\,\,

		\item[\rm\bf ii)]{\sc (adjoint differential inclusion)} $$\frac{dp}{ds}\in - \partial_x^C H\left(\ol{y},p_0,p,\lambda,\pi,\ol{w}^0,\ol{w},\ol{\alpha}
		\right);$$
		\item[\rm\bf iii)]{\sc(non transversality)} $$\ds(p_0,p(\ol{S}))\in -\lambda \partial_{(t,x)}^C\Psi\Big((\ol{y}^0,\ol{y})(\ol{S})\Big)-\ol{\bigcup_{\mathcal{T}\in \mathbcal{T}}\mathcal{T}^\perp}\,\,;$$

		\item[\rm\bf iv)]{\sc (first order maximization)}\,\,\, For almost all $s \in [0,\ol{S}]$,
		$$\max_{(w^0,w,a) \in \xR_+\times C \times A}\Big[H(\ol{y}(s),p_0,p(s),\lambda,0,w^0,w,a)\Big]=H(\ol{y}(s),p_0,p(s),\lambda,0,\ol{w}^0(s),\ol{w}(s),\ol{\alpha}(s)).
		$$
		
		\item[\rm\bf v)]{\sc(nonsmooth Goh condition)} \,\,\, If, in addition,  $\|\ol{w}\|_1<K \text{ and} \ \hat{l}_1( \cdot,0)\equiv 0$, then
		\bel{thlie}
		{0\in p(s)\, [g_i,g_j]_{set}(\ol{y}(s))  \qquad i,j \in \{1,\ldots,m_1\}, \quad \text{for a.e.} \ t \in [0,\ol{T}]. }
		\eeq
		
	\end{itemize}
\end{thrm}

The proof of  Theorem \ref{TeoremaPrincipale} will be given in Section \ref{proofsec}. 
%It includes the use of Theorem \ref{qdqth}, whose   proof   is given in subsection \ref{th3sec}.

\begin{rmrk}\label{iiisimplified}{  Let us point out that if the  {Clarke's  \it    tangent cone} to the target  happens to be a $QDQ$-approximating cone as well, the set of conditions   i)-iv) coincides with a (first order) non-smooth Pontryaging Maximum Principle of the kind one finds in several  books (see e.g. \cite{Vinter} e\cite{Clarke2} )} The same of course can be said for the smooth case where the target is a differential submanifold, in which case the tangent space to the target is automatically a $QDQ$-cone. Furthermore, in the smooth case,\eqref{thlie} coincides with the classical Goh condition.
	
\end{rmrk}

%{\fra So che cia hai lavorato, ma ho tolto il caso del minimo non impulsivo, perché questo era stato aggiunto per trattare più facilmente l' esempio. Ma l' esempio adesso è tutto fatto con il pseudotempo $s$, come nel teorema principale, e in fondo è meglio così, per non aggiungere inutile carne al fuoco }

\section{Proof of  the Maximum Principle}\label{proofsec}
\subsection{An equivalent fixed end-time problem}

For every process $({S},{w}^0,{w},{\alpha},{y}^0,{y},{\beta})$, we will set  $$\ol{y}^l(s):=\displaystyle\int_0^sl^e(\overline{y}(\sigma),\ol{w}^0(\sigma),\ol{w}(\sigma),\ol{\alpha}(\sigma)) d\sigma$$  (so that  $\ol{y}^l$  is  the unique Carathéodory solution to the trivial differential equation $\ds\frac{dy^l}{ds}(s)=l^e(\overline{y}(s),\ol{w}^0(s),\ol{w}(s),\ol{\alpha}(s))$ with initial condition $y^l(0)=0$.)\\

Let us begin with a further (and standard) reparametrization procedure which allows us to reduce  problem  $(P_{ext})$ to a problem with a fixed end-time. 

Let us fix  $\bar S>0$, $\rho>0$. We say that $(\bar S,w^0,w,\alpha, \zeta, y^0,y,y^l,\beta)$ is a {\it  rescaled space-time process}  if  $$(\bar S,w^0,w,\alpha,\zeta)(\cdot) \in \mathcal{W} \times L^\infty([0,\ol{S}],[-\rho,\rho])$$ and $\big((y^0,y,y^l),\beta\big)$ is the unique (Carathéodory) solution of {\it the rescaled Cauchy problem }

\begin{equation} \label{sistema}
	\begin{cases}
		\ds\frac{d}{ds}\big((y^0,y,y^l),\beta\big)=\Big(\cl{F}(y,w^0,w,a),|w|\Big) \cdot \Big(1+\zeta\Big)\qquad s\in[0,\bar S]\\
		\big((y^0,y,y^l),\beta\big)(0)=\big((0,\hat{x},0),0\big)
	\end{cases}
\end{equation}
Moreover, we call $(\bar S,w^0,w,\alpha, \zeta, y^0,y,y^l,\beta)$ {\it  feasible} as soon as   $\big((y^0,y),\beta\big) \in \mathfrak{T}\times [0,K]$.\\
The {\it  rescaled optimization problem} is defined as  \begin{equation}\label{problemariscalato}\begin{cases}
		\min \left( \Psi((y^0,y)(\ol{S}))+y^l(\ol{S})\right),\\
		\text{ over feasible rescaled processes}.
\end{cases}\end{equation}
With standard arguments one shows that,  for small enough $\rho>0$, a  process  $\opprocess$ is a %canonical
 local minimizer for the extended  problem  $(P_{ext})$ if and only if the  rescaled space-time process $(\bar S,\bar w^0,\bar w, \bar\alpha, 0, \bar y^0,\bar y,\bar y^l,\bar \beta)$  is a local minimizer for fixed-end-time problem \eqref{problemariscalato}.  \footnote{ Actually, the role of the auxiliary parameter $\zeta$ is fictitious, because of the rate-independence of problem $(P_{ext})$.  However we use it, for it makes proofs simpler.}  

Therefore, in the proof of the Maximum Problem we are allowed  to  {\it replace the hypothesis of the theorem with the following one:}

\begin{itemize}
	\item {\it The process $(\bar S,\bar w^0,\bar w, \bar\alpha, \bar \zeta\equiv 0 , \bar y^0,\bar y,\bar y^l,\bar \beta)$ is a local minimizer of the rescaled problem \eqref{problemariscalato}. \it}
	
\end{itemize}

\subsection{Set separation}
For some $\delta>0$, let us  consider the {\it $\delta$-reachable set}
$$\mathfrak{R}_{\delta}:=\left\lbrace \begin{aligned}& \Big(y^0,y,y^l +\Psi (y^0,y),\beta\Big)(\ol{S}):\,\,\, (\ol{S},w^{0},w,\alpha,\zeta,y^0,y,y^l ,\beta) \text{ is a rescaled }\\ &\text{process that  verifies }  \|(y^0-\ol{y}^0,y-\ol{y},y^l-\ol{y}^l,\beta-\ol{\beta})\|_{\infty} < \delta \end{aligned}\right \rbrace\subseteq \xR^{1+n+1+1}    $$
and the {\it projected $\delta$-reachable set }
$$\mathfrak{R}^{'}_{\delta}:= \prj\Big(\mathfrak{R}_{\delta}\Big)\subseteq \ \xR^{1+n+1} , $$
where the projection operator $\prj$ is defined by setting $\prj (x^0,x,x^l,\beta):=  (x^0,x,x^l),$ for all $  (x^0,x,x^l,\beta)\in \xR^{1+n+1+1}.$
Let us introduce also the {\it profitable set}
$$
\mathfrak{P}:=\Bigg(\Big(\mathfrak{T} \times \Big]-\infty, \ol{y}^l(\ol{S}) +\ol{\Psi}(\ol{S}) \Big[\,\Big)\bigcup \Big\{\big(\ol{y}^0, \ol{y}, \ol{y}^l(\ol{S}) +\ol{\Psi}(\ol{S})\big)\Big\}\Bigg)\times[0,K]
$$
where $\overline{\Psi}(s):=\Psi(\ol{y}^0(s),\ol{y}^{ }(s))$ for all $s$, and the {\bf projected profitable set}
$$
\qquad\qquad\mathfrak{P}^{'} :  =  \prj\Big(\mathfrak{P}\Big) = \Bigg(\mathfrak{T} \times\Big]-\infty, \ol{y}^l(\ol{S}) +\ol{\Psi}(\ol{S}) \Big[\Bigg)\bigcup \Big\{\big(\ol{y}^0, \ol{y}, \ol{y}^l(\ol{S}) +\ol{\Psi}(\ol{S})\big)\Big\}
$$

\begin{lmm}\label{setsep} 	Let $\opprocess$  as in Theorem \ref{TeoremaPrincipale}, and  let us assume that $\ol{\beta}(\ol{S})<K$. Then for any  $\delta>0$ sufficiently small, the projected profitable set ${\mathfrak{P}^{'}}$ and the projected $\delta$-reachable set  $ \mathfrak{R}^{'}_{\delta}$ are locally separated at $\big(\ol{y}^0, \ol{y}, \ol{y}^l(\ol{S}) +\ol{\Psi}(\ol{S})\big)$.
	
\end{lmm}  
\begin{proof}Indeed, by the definition of local minimizer it follows that  the  profitable set ${\mathfrak{P}}$ and the $\delta$-reachable set $\mathfrak{R}_{\delta}$ are locally separated at  $\big(\ol{y}^0, \ol{y}, \ol{y}^l(\ol{S}) +\ol{\Psi}(\ol{S}), \ol{\beta}(\ol{S})\big)$. From this one gets the thesis trivially.
	%(see  \cite{ArMotRamp}, Lemma 6.12).
\end{proof}

\subsection{Finitely many variations}

With the ultimate aim of applying a suitable  separability criterion for approximating cones, we  now build  a $QDQ$-approximating multicone to the  projected $\delta$-reachable set  $ \mathfrak{R}^{'}_{\delta}$ at  $\big(\ol{y}^0, \ol{y}, \ol{y}^l(\ol{S}) +\ol{\Psi}(\ol{S})\big)$.
Let us define the set $\mathfrak{V}$ of {\it  variation generators} as the union
$\mathfrak{V}:=\mathfrak{V}_{ndl}\bigcup \mathfrak{V}_{brk}$,
where $\mathfrak{V}_{ndl}$ and $ \mathfrak{V}_{brk}$ are the sets
of {\it  needle variation generators} and of {\it  bracket-like variation generators} defined as
$\mathfrak{V}_{ndl}:=\xR_+ \times C \times A \times [-\rho,\rho]$ and  $\mathfrak{V}_{brk} :=\Biggl[\{1,\ldots,m_1\}^2\setminus \text{diag}\Big(\{1,\ldots,m_1\}^2\Big)\Biggr],$  respectively.
%{\fra anche qui c'e' $W$ non definito}
\begin{dfntn}Let $(0,\ol{S})_{Leb}\subset [0,\ol{S}]$ be the set of Lebesgue points of the function $s \mapsto (\ol{w}^0(s),\ol{F}^e(s),\ol{l}^e(s),|\ol{w}|(s))$,  where $\ol{F}^e$ and $\ol{l}^e$ denote the functions $F^e$ and $l^e$ evaluated along the optimal process $(\bar S,\bar w^0,\bar w, \bar\alpha, \bar \zeta\equiv 0 , \bar y^0,\bar y,\bar y^l,\bar \beta)$  of the rescaled problem. For every variation generator $\mathbf{c}\in \mathfrak{V}$, let us define the variation vector  
	$\left(v^0_{\mathbf{c},\ol{s}}, v^{ }_{\mathbf{c},\ol{s}}, v^l_{\mathbf{c},\ol{s}} \right)$  at an instant $\ol{s}$ as follows:\footnote{As in the standard maximum principle, the fact of not considering   pairs $(\mathbf{c},\ol{s}) \in 
		\mathfrak{V}_{ndl}\times\Big((0,\ol{S})\backslash (0,\ol{S})_{Leb}\Big) $ is completely irrelevant, in that $(0,\ol{S})\backslash (0,\ol{S})_{Leb})$ has zero measure. }	
	$$
	\left(v^0_{\mathbf{c},\ol{s}}, v^{ }_{\mathbf{c},\ol{s}}, v^l_{\mathbf{c},\ol{s}} \right):= \left\{
	\begin{array}{ll}\begin{array}{l}
			\left\{\begin{pmatrix} w^0(1+\zeta)-\ol{w}^0(\ol{s})\\ F^e(\ol{y}(\ol{s}),w^0,w,a)(1+\zeta)-\ol{F}^e(\ol{s})\\ l^e(\ol{y}(\ol{s}),w^0,w,a)(1+\zeta)-\ol{l}^e(\ol{s})  \end{pmatrix}\right\}
		\end{array}\,\,\,  &\begin{array}{l}\text{if}\,\,
			\mathbf{c}=(w^0,w,a,\zeta)\in \mathfrak{V}_{ndl}, \\\,\,\,\,\, \text{and}\,\,\, \ol{s} \in (0,\ol{S})_{Leb}\end{array}
		\\[12mm]
		\{0\} \times [g_i,g_j]_{set}(\ol{y}(\ol{s}))\times\{0\}\qquad\qquad &\begin{array}{l}\,\,\, \text{if}
			\,\, \mathbf{c}=(i,j)\in \mathfrak{V}_{brk} \\\,\,\,\,\,\text{and}\,\,\, \ol{s} \in (0,\ol{S}).
		\end{array}
	\end{array}\right.
	$$
	Moreover, when $
	\mathbf{c}=(w^0,w,a,\zeta)\in \mathfrak{V}_{ndl}, $ and $ \ol{s} \in (0,\ol{S})_{Leb} $, we set
	$
	v^\nu_{\mathbf{c},\ol{s}} := |w|(1+\zeta)-|\ol{w}(\ol{s})|.
	$

\end{dfntn}
Let us point out that, to retain uniformity of notation, we always regard $\left(v^0_{\mathbf{c},\ol{s}}, v_{\mathbf{c},\ol{s}}^{ }, v^l_{\mathbf{c},\ol{s}} \right)$ as a {\it subset of vectors of  $\xR^{1+n+1}$ }, though, as soon as $\mathbf{c}\in \mathfrak{V}_{ndl}$,  it reduces to the singleton formed by the usual needle variation vector.

\begin{dfntn}
	Let us fix a rescaled control  ${\mathbf{w}}=({w}^0,{w},{\alpha},{\zeta}) \in L^\infty\big([0,\ol{S}],\xR_+\times C\times A\times [-\rho,\rho]\big)$ (with  $\essinf(w^0+|w|)>0$) and an instant   $\ol{s} \in (0,\ol{S})$.

	\begin{itemize}
		\item If  $\mathbf{c}=(\hat w^0,\hat w,\hat a,\hat\zeta)\in \mathfrak{V}_{ndl}$, we call  {\it  needle control variation of ${\mathbf{w}}$ at $\ol{s}$ associated to $\mathbf{c}$} the family  of controls $\big\{ {\mathbf{w}}_{\varepsilon,\mathbf{c},\ol{s}}(s): \, \varepsilon \in [0,\ol{s})\big\}$ defined as  $${\mathbf{w}}_{\varepsilon,\mathbf{c},\ol{s}}(s)=\begin{cases}{\mathbf{w}}(s) &\text{ if } s \in [0,\ol{s}-\varepsilon) \cup (\ol{s},\ol{S}]\\ (\hat w^0,\hat w,\hat a,\hat\zeta) &\text{ if } s \in [\ol{s}-\varepsilon,\ol{s}].\end{cases}$$   
		
		\item 
		If 
		$\mathbf{c}=(i,j)\in  \mathfrak{V}_{brk}$, we call  {\it bracket-like variation  of $\mathbf{{w}}$ at $\bar s$} the family  $\Big\{ {\mathbf{w}}_{\varepsilon,\mathbf{c},\ol{s}}(s): \, 0<8\sqrt{\varepsilon}\le\ol{s}\Big\}$ of controls defined as
		$$\mathbf{{w}}_{\varepsilon, \mathbf{c},\ol{s}}(s)=\begin{cases}
			\mathbf{{w}}(s) & \text{ if } s \not \in [\ol{s}-8\sqrt{\varepsilon},\ol{s}]\\
			(2{w}^0,2{w},{\alpha},{\zeta})\circ \gamma^\varepsilon(s) & \text{ if }s \in [\ol{s}-8\sqrt{\varepsilon},\ol{s}-4\sqrt{\varepsilon}]\\
			(0,\mathbf{e}_i,a,0)& \text{ if }s \in [\ol{s}-4\sqrt{\varepsilon},\ol{s}-3\sqrt{\varepsilon}]\\
			(0,\mathbf{e}_j,a,0)& \text{ if }s \in [\ol{s}-3\sqrt{\varepsilon},\ol{s}-2\sqrt{\varepsilon}]\\
			(0,-\mathbf{e}_i,a,0)& \text{ if }s \in [\ol{s}-2\sqrt{\varepsilon},\ol{s}-\sqrt{\varepsilon}]\\
			(0,-\mathbf{e}_j,a,0)& \text{ if }s \in [\ol{s}-\sqrt{\varepsilon},\ol{s}],
		\end{cases}$$
		where $a\in A$ is arbitrarily chosen\footnote{Since $w^0=0$, the choice of $a$ is indeed irrelevant.} 
		and 
		$\gamma^\varepsilon(s):=2s-\ol{s}+8\sqrt{\varepsilon}$, 
		%$\forall s \in  [\ol{s}-8\sqrt{\varepsilon},\ol{s}-4\sqrt{\varepsilon}],$ 
		
	\end{itemize}
\end{dfntn}
\subsection{ $QDQ$-approximating cones to $ \mathfrak{R}^{'}_{\delta}$} \label{th3sec}
For any $(\ol{s},\mathbf{c})\in [0,S] \times \mathfrak{V}$ and any $\varepsilon$ sufficiently small, consider the functional  $\mathcal{A}_{\varepsilon,\mathbf{c},\ol{s}}$ (from the space of rescaled controls ${\mathbf{w}}$ into itself) defined by setting  $\mathcal{A}_{\varepsilon,\mathbf{c},\ol{s}}(\mathbf{w}):=\mathbf{w}_{\varepsilon,\mathbf{c},\ol{s}}.$
In addition, given $N$ variation generators $\mathbf{c}_1,\ldots,\mathbf{c}_N \in \mathfrak{V}$ and $N$ instants $ 0<{s}_{1} < {s}_2 < \ldots {s}_N\leq \ol{S}$ for a $\tilde\varepsilon>0$ sufficiently small, let us define the multiple variation  $$[0,\tilde\varepsilon]^N\ni \small{\bm{\varepsilon}}\mapsto \mathbf{\ol{w}}_{\bm{\varepsilon}}:=\mathcal{A}_{\varepsilon_N,\mathbf{c}_N,{s}_N} \circ \ldots \circ \mathcal{A}_{{\varepsilon_1},\mathbf{c}_1,{s}_1} (\mathbf{\ol{w}}).$$
Let us set  $(\ol{w}^0_{\bm{\varepsilon}},\ol{w}_{\bm{\varepsilon}},\ol{a}_{\bm{\varepsilon}},\ol{\zeta}_{\bm{\varepsilon}}):=\mathbf{\ol{w}}_{\bm{{\bm{\varepsilon}}}}$, and let  us use  $ (y^0_{\bm{\varepsilon}},y_{\bm{\varepsilon}},y^l_{\bm{\varepsilon}},\beta_{\bm{\varepsilon}})  $ to denote  the solution (on $[0,\ol{S}]$) of the Cauchy problem \footnote{Of course, $\mathbf{\ol{w}}_{\bm{\varepsilon}}$ and $(\extpbf{y},\beta_{\bm{\varepsilon}} )$ depend also on the parameters $\mathbf{c}_k$ and ${s}_k$, but we avoid writing them when possible in order to simplify the notation.}
\begin{equation}
\label{multiperturbato} \tag{$P_{\varepsilon,\mathbf{c},\ol{s}}$}
\begin{cases}
\ds\frac{d}{ds} (y^0,y,y^l,\beta)=\Bigg(\cl{F}\Big(y,\ol{w}^0_{\bm{\varepsilon}}, \ol{w}_{\bm{\varepsilon}},\ol{a}_{\bm{\varepsilon}}\Big),|w_{\bm{\varepsilon}}| \Bigg) (1+\ol{\zeta}_{\bm{\varepsilon}})\\[3mm]
(y^0,y,y^l,\beta)(0)=(0,\hat{x},0,0)
\end{cases}
\end{equation} \\ 
\begin{lmm} \label{lemmapseudoaffine}
	The map $\mathbf{Y}:\xR_+^N\to \xR^q$ defined by setting
	$$\mathbf{Y}(\bm{\varepsilon})  :=
	\Big(y^0_{\bm{\varepsilon}}(\ol{S})\,,\,
	y_{\bm{\varepsilon}}(\ol{S})\,,\,\,y^l_{\bm{\varepsilon}}(\ol{S}) \Big)
	$$   
	satisfies  the hypothesis \eqref{pseudoaffine} with $F=\mathbf{Y}$ and  $q:= 1+n+1$,
	namely,
	% there exists a function $g(\bm{\varepsilon}) = o(|\bm{\varepsilon}|)$ 
	one has
	\begin{equation}\label{pseudoaffine2}
		\mathbf{Y}(\bm{\varepsilon})-\mathbf{Y}(0)=\sum\limits_{i=1}^{N} \big(\mathbf{Y}(\varepsilon_i \mathbf{e}_i)- \mathbf{Y}(0)\big)+ o(|\bm{\varepsilon}|),\qquad \forall \bm{\varepsilon} =(\varepsilon^1,\dots,\varepsilon^N)\in\xR_+^N.\end{equation}
\end{lmm}
\begin{proof}
	Let $\eta: \xR^n \to \xR_+$  and, for every  $\delta>0$, $\eta_\delta$  be a $C^\infty$ mollifier as in Proposition \ref{Clarketh}. For any  $\delta>0$,
	let us define the mollified vector field  $$\cl{F}_{\delta}\Big(y,w^0,w,a\Big):=\int_{\xR^n}\cl{F}\Big(y+h,w^0,w,a\Big) \eta_{\delta}(h)\,dh.$$ Observe that  the control vector field  $\cl{F}$
	is  continuous, and, in addition, it is locally Lipschitz in the variable $y$. Moreover, we can apply a cut off technique  and make $\cl{F}$ and $\Psi$ equal to zero outside  a compact set containing a small neighbourhood of our local minimizer, so that {\it we can assume that  $\cl{F}$ and $\Psi$ are  globally Lipschitz as well.} It follows that $\cl{F}_{\delta}$ converges uniformly to $\cl{F}$ as $\delta$ goes to $0$.
	For any  fixed $\bm{\varepsilon}\in \xR_+$ with a   suitably small norm, let us introduce the mollified  {Cauchy} problem 
	
	\begin{equation} \label{mollificatoperturbato}
		\begin{cases}
			\ds\frac{d}{ds}\big((y^0,y,y^l),\beta\big)=\Big(\cl{F}_\delta(y,w^0_{\bm{\varepsilon}},w_{\bm{\varepsilon}},a_{\bm{\varepsilon}}),|w|\Big) \cdot \Big(1+\zeta_{\bm{\varepsilon}}\Big)\\\\
			\big((y^0,y,y^l),\beta\big)(0)=\big((0,\hat{x},0),0\big)
		\end{cases}
	\end{equation}
	and let us use $ \left(y^0_{\delta,\bm{\varepsilon}},y_{\delta,\bm{\varepsilon}},y^l_{\delta,\bm{\varepsilon}},\beta_{\bm{\varepsilon}}\right)$ to denote its  unique solution.\\
	We also set $$\mathbf{Y}_{\delta}(\bm \varepsilon)	:=
	\Big(y^0_{\delta,\bm{\varepsilon}}(\ol{S})\,\,,\,\,
	y_{\delta,\bm{\varepsilon}}(\ol{S})\,\,,\,\,y^l_{\delta,\bm{\varepsilon}}(\ol{S}) + \Psi\left(y^0_{\delta,\bm{\varepsilon}}(\ol{S}),
	y_{\delta,\bm{\varepsilon}}(\ol{S})\right) \Big)$$
	
	Let us define the function $z_{\delta,\bm{\varepsilon}}(s):=\left|(y^0_{\delta,\bm{\varepsilon}},y_{\delta,\bm{\varepsilon}},y^l_{\delta,\bm{\varepsilon}})(s)-(y^0_{\bm{\varepsilon}},y_{\bm{\varepsilon}},y^l_{\bm{\varepsilon}})(s)\right|$, $s\in [0,\bar S]$, and let us observe that, from the inequality  $$\begin{array}{c}
		\ds z_{\delta,\bm{\varepsilon}}(s)\le \int_0^s \Big|\cl{F}_{\delta}(y_{\delta,\bm{\varepsilon}},w^0_{\bm{\varepsilon}},w_{\bm{\varepsilon}}
		,a_{\bm{\varepsilon}})-\cl{F}(y_{\bm{\varepsilon}},
		w^0_{\bm{\varepsilon}},w_{\bm{\varepsilon}},a_{\bm{\varepsilon}})\Big|\Big(1+\zeta_{\bm{\varepsilon}}\Big)d\sigma \le \\\ds  \int_0^s \Big|\cl{F}_{\delta}(y_{\delta,\bm{\varepsilon}},w^0_{\bm{\varepsilon}},w_{\bm{\varepsilon}}
		,a_{\bm{\varepsilon}})-\cl{F}(y_{\delta,\bm{\varepsilon}},w^0_{\bm{\varepsilon}},w_{\bm{\varepsilon}}
		,a_{\bm{\varepsilon}})\Big|\Big(1+\zeta_{\bm{\varepsilon}}\Big)d\sigma +L(1+2\rho)\int_0^s z_{\delta,\bm{\varepsilon}}(\sigma) \le\\ 2K(1+2\rho)\overline{S}\delta+L(1+2\rho)\int_0^s z_{\delta,\bm{\varepsilon}}(\sigma)d\sigma,\quad 
	\end{array} \quad \footnote{ $K$ is a bound for the maps $\cl{F}$ and $\cl{F}_\delta$  and $L$ is a Lipschitz constant for the maps  $(y^0,y,y^l)\mapsto\cl{F}(y^0,y,y^l,w^0,w
		,a)$, independent of $(w^0,w,a)$}$$
	and Gronwall's Lemma, we deduce\begin{equation}\label{convergenza} \left|(y^0_{\delta,\bm{\varepsilon}},y_{\delta,\bm{\varepsilon}},y^l_{\delta,\bm{\varepsilon}})(s)-(y^0_{\bm{\varepsilon}},y_{\bm{\varepsilon}},y^l_{\bm{\varepsilon}})(s)\right|=  z_{\delta,\bm{\varepsilon}}(s)\le C\delta,\quad \forall s\in [0,\bar S],\end{equation}  where $C$ is a positive constant depending only on $\overline{S}$, $K$, and $L$. By choosing $\delta=\delta(|\bm\varepsilon|)=|\bm\varepsilon|^2$, we get  $\mathbf{Y}(\bm \varepsilon)=\mathbf{Y}_{|\varepsilon^2|}(\bm\varepsilon)+o(|\bm{\varepsilon}|).$ 
	%Now, since  systems \eqref{mollificatoperturbato} are smooth, one can elementarily prove (see e.g. \cite{ArMotRamp}, Lemma 5.4) that, for any positive $\delta$ sufficiently small,  the function $\mathbf{Y}_\delta$ verifies condition 
%	\eqref{pseudoaffine2}, i.e.
%	$$
%	\mathbf{Y}_\delta(\bm{\varepsilon})-\mathbf{Y}_\delta(0)=\sum\limits_{i=1}^{N} \big(\mathbf{Y}_\delta(\varepsilon_i \mathbf{e}_i)- \mathbf{Y}_\delta(0)\big)+ g_\delta({\varepsilon}),\qquad \forall \bm{\varepsilon} =(\varepsilon^1,\dots,\varepsilon^N)\in\xR^N,
%	$$ with  $g_\delta(\bm{\varepsilon})= o(|\bm{\varepsilon}|).$ {\fra QUESTO ERA RIMASTO DA COSE VECCHIE? 
	We have reduced to a smooth system, hence the fact that a function like our $\mathbf{Y}_{|\varepsilon^2|}$ has the desired property \eqref{pseudoaffine2} is something very well known in the literature. The thesis follows automatically as $\mathbf{Y}_{|\varepsilon^2|}(\varepsilon_i \mathbf{e}_i)$ is again distant at most $o(|\bm{\varepsilon}|)$ from $\mathbf{Y}(\varepsilon_i \mathbf{e}_i)$.
	As a matter of fact, 
	$$\begin{array}{c}	\mathbf{Y}(\bm{\varepsilon})-\mathbf{Y}(0)= 
		\big(\mathbf{Y}(\bm{\varepsilon})-\mathbf{Y}_{|\bm{\varepsilon}|^2}(\bm{\varepsilon})\big) + \big(\mathbf{Y}_{|\bm{\varepsilon}|^2}(\bm{\varepsilon})-\mathbf{Y}_{|\bm{\varepsilon}|^2}(0)\big) +
		\big(\mathbf{Y}_{|\bm{\varepsilon}|^2}(0) - \mathbf{Y}(0)\big) =\\\
		\ds	\sum\limits_{i=1}^{N} \big(\mathbf{Y}_{|\bm{\varepsilon}|^2}(\varepsilon_i \mathbf{e}_i)- \mathbf{Y}_{|\bm{\varepsilon}|^2}(0)\big)+\Big(o(|\bm{\varepsilon}|)	+ 	\big(\mathbf{Y}(\bm{\varepsilon})-\mathbf{Y}_{|\bm{\varepsilon}|^2}(\bm{\varepsilon})\big)  +
		\big(\mathbf{Y}_{|\bm{\varepsilon}|^2}(0) - \mathbf{Y}(0)\big)\Big) =  \\
		\ds	\sum\limits_{i=1}^{N} \big(\big(\mathbf{Y}(\varepsilon_i \mathbf{e}_i)- \mathbf{Y}(0)\big)+o(|\bm{\varepsilon}|)
	\end{array}$$

\end{proof}

\begin{dfntn}
	Let  $N$ be  a natural number,  and  let us choose $N$ variation generators $\mathbf{c}_1,\ldots,\mathbf{c}_N \in \mathfrak{V}$ and $N$ instants $0<{s}_{1} < {s}_2 < \ldots \le {s}_N\leq \ol{S}$, with $s_k\in  [0,\ol{S}]_{Leb} $ as soon as $\mathbf{c}_k \in \mathfrak{V}_{ndl}$. 
	For every $k=1,\ldots,N$, any $L^1$-map $[0,\ol{S}]\ni s\mapsto ({M},\omega)(s)\in \Lin(\xR^n,\xR^n)\times(\xR^n)^*$, and any $(m_t,m_x) \in (\mathbb{R}^{n+1})^*$ let us  consider the $(1+n+1)\times(1+n+1)$  matrix
	%\begin{equation}\label{matrix}
	
	$$
	\mathcal{E}'_k(m_t,m_x,{M},\omega) :=\begin{pmatrix} 1&0&0\\ 0_{n\times 1}  &e^{\bigintssss_{{s}_k}^{\ol{S}}{M}\,} & 0_{n\times 1} \\ m_t & m_x e^{\bigintssss_{{s}_k}^{\ol{S}}{M}\,}+\displaystyle{\bigintssss_{{s}_k}^{\ol{S}}\Big(\omega(s)  e^{\bigintssss_{s}^{{s}_k}{M}}\,}\Big)ds \,\,\,\,\,\,\,\,&1\end{pmatrix}
	$$
	{ (which transports  vectors from the tangent space at $(\ol{y}^0,\ol{y},\ol{y}^l)(s_k)$ to the  tangent space at $(\ol{y}^0,\ol{y},\ol{y}^l)(\ol{S})$)} and,
	in the special case when $\bm{c}_k \in \mathfrak{V}_{ndl}$ $\forall k\in \{1,\ldots,N\}$,
	the $(1+n+1+1)\times(1+n+1+1)$ matrix %\begin{equation}\label{matrixgrande}
	$$
	\mathcal{E}_k(m_t,m_x,{M},\omega) :=\begin{pmatrix}1  &  0   & 0 & 0 \\ 0_{n\times 1}   &e^{\bigintssss_{{s}_k}^{\ol{S}}{M}\,} &0_{n\times 1}  & 0\\ m_t &m_x e^{\bigintssss_{{s}_k}^{\ol{S}}{M}\,}+\displaystyle{\bigintssss_{{s}_k}^{\ol{S}}\Big(\omega(s)  e^{\bigintssss_{s}^{{s}_k}{M}}\,}\Big)ds \,\,\,\,\,\,\,\,&1 & 0 \\ 0&0&0&1\end{pmatrix}
	%\end{equation}
	$$

	where the exponential of a matrix is defined as in Subsection \ref{notations} and $0_{n\times 1}$ stands for a column of $n$ zeros.  Subsequently,  let the subsets $\Lambda'_N \subset  \Lin(\xR^N,\xR^{1+n+1})$  $\Lambda_N \subset  \Lin(\xR^N,\xR^{1+n+1+1})$ be  defined as
	$$ \Lambda'_N \!:= \!\left\lbrace 
	\left(\mathcal{E}'_{1} (m_t,m_x,{M,\omega}) \begin{pmatrix}V_1^0\\[3mm] V_1 \\[3mm]V_1^l\end{pmatrix},\dots,\mathcal{E}'_{N}(m_t,m_x,{M,\omega})\begin{pmatrix}V_N^0
		\\[3mm] V_N \\[3mm]V_N^l\end{pmatrix}\right)  , \quad 
	\begin{array}{l} ({M},\omega)(\cdot)\,  \text{ is a measurable}\\ \text{selection of}  
		\,\,	\partial_{y}\big(\ol{F}^e,\ol{l}^e\big)(\cdot ), \\ (m_t,m_x) \in \partial_{(t,x)}^C\Psi((\ol{y}^0,\ol{y})(\ol{S}))\\  
		(V_k^0, V_k,V_k^l) \in (v_{\mathbf{c}_k,{s}_k}^0, v_{\mathbf{c}_k,{s}_k}, v^l_{\mathbf{c}_k,{s}_k})\\ \forall k=1,\ldots,N \end{array}
	\right \rbrace,$$
	$$ \Lambda_N \!:= \!\left\lbrace 
	\left(\mathcal{E}_{1} (m_t,m_x,{M,\omega}) \begin{pmatrix}V_1^0\\[3mm] V_1 \\[3mm]V_1^l \\[3mm]V_1^\nu\end{pmatrix},\dots,\mathcal{E}_{N}(m_t,m_x,{M,\omega})\begin{pmatrix}V_N^0
		\\[3mm] V_N \\[3mm]V_N^l\\[3mm]V_N^\nu\end{pmatrix}\right)  , \,
	\begin{array}{l} \text{where}\,\, ({M},\omega)(\cdot), m_t, m_x \\\text{ and} 
		(V_k^0, V_k,V_k^l)\\ \text{are as in the} \\ \text{definition of }\Lambda'_N\end{array}
	\right \rbrace.$$
	
\end{dfntn}

Corollary  \ref{qdqcor} below represents the most important technical step  of the proof of our  maximum principle. It is a straightforward consequence of the following result:

\begin{thrm}\label{qdqth}
	Let  $\big(\ol{y}^0, \ol{y}, \ol{y}^l,\ol{\beta}\big)$ and $\big({y}^0_{\bm{\varepsilon}}, {y}_{\bm{\varepsilon}}, {y}^l_{\bm{\varepsilon}},{\beta}_{\bm{\varepsilon}}\big)$ as above. If we assume the  extra assumption $\hat{l}_1( \cdot,0,\cdot )\equiv 0$,
	%Let $N$ be any natural number and consider any $\mathbf{c}_1,\ldots,\mathbf{c}_N \in \mathfrak{V}$ and suitable times ${s}_1<{s}_2<\ldots<{s}_N \in[0,\ol{S}]$.\\
	then the set  ${\Lambda}'_N$ is a $QDQ$ at $\bm{0}$ {of the map $$\mathbf{Z}:=
		\Big(y^0_{\bm{\varepsilon}}(\ol{S})\,,\,
		y_{\bm{\varepsilon}}(\ol{S})\,,\,\,y^l_{\bm{\varepsilon}}(\ol{S}) + \Psi\left(y^0_{\bm{\varepsilon}}(\ol{S}),
		y_{\bm{\varepsilon}}(\ol{S})\right) \Big)
		$$}
	in the direction of $\xR_+^{N}$.

	Moreover, in the special case when $\mathbf{c}_k \in \mathfrak{V}_{ndl}$ for all $k\in \{1,\ldots,N\}$  (and $\hat{l}_1( \cdot,0, )$ is possibly non vanishing), $\Lambda_N$ is a $QDQ$ at $\bm{0}$ of the map $\xR_+^N\ni\bm{\varepsilon}  \mapsto \left(\mathbf{Z}(\bm{\varepsilon}), \beta_{\bm{\varepsilon}}\right)$ (in the direction of $\xR_+^{N}$).

\end{thrm} 

\begin{crllr}\label{qdqcor} Let us use the same notations as in Theorem \ref{qdqth} and let us assume that $\hat{l}_1( \cdot,0, )\equiv 0.$ 
	For any choice of $\delta>0$,	the set $$\Lambda'_N \xR_+^N:=\Big\{{L}'  \xR_+^N:\, L' \in \Lambda'_N\Big\}$$   is a $QDQ$-approximating multicone of  the projected $\delta$-reachable set $\mathcal{R}'_\delta$ at $\big(\ol{y}^0, \ol{y}, \ol{y}^l(\ol{S}) +\ol{\Psi}(\ol{S})\big)$. Moreover, in the special case when $\mathbf{c}_k \in \mathfrak{V}_{ndl}$ for all $k\in \{1,\ldots,N\}$ (and $\hat{l}_1( \cdot,0, )$ is possibly non vanishing), the set $$\Lambda_N \xR_+^N:=\Big\{L  \xR_+^N:\, L \in \Lambda_N\Big\}$$  is a $QDQ$-approximating multicone of the $\delta$-reachable set $\mathcal{R}_\delta$ at $\big(\ol{y}^0, \ol{y}, \ol{y}^l(\ol{S}) +\ol{\Psi}(\ol{S}),\ol{\beta}(\ol{S})\big)$.
\end{crllr}

\begin{proof}[Proof of Theorem \ref{qdqth} ] {First of all, notice that $\mathbf{Z}$ is obtained by composition of the map $\mathbf{Y}$ we defined in Lemma \ref{lemmapseudoaffine}   with the function $(t,x,c)\mapsto (t,x,c+\Psi(t,x))$. This means that in order to build a Quasi Differential Quotient of $\mathbf{Z}$ it is enough to build a $QDQ$ of $\mathbf{Y}$. Indeed,  thanks to Proposition \ref{Clarketh}, the Clarke's Generalized Jacobian of $(t,x,c)\mapsto (t,x,c+\Psi(t,x))$ is a also  $QDQ$ (in the direction of $\xR^{1+n+1}$) of  $(t,x,c)\mapsto (t,x,c+\Psi(t,x))$. In turn, the Clarke's Generalized Jacobian of $(t,x,c)\mapsto (t,x,c+\Psi(t,x))$ at $(\ol{S},\ol{y}(\ol{S}),\ol{y}^l(\ol{S}))$ is easily seen to coincide with  $$\left\{\begin{pmatrix}1 & \mathbf{0} &0 \\ \mathbf{0} & \mathbf{1} & \mathbf{0} \\ m_t & m_x & 1 \end{pmatrix}, (m_t,m_x) \in \partial_{(t,x)}^C \Psi(\ol{S},\ol{y}(\ol{S}))\right\}$$
		(where $\partial_{(t,x)}^C \Psi$ is the Clarke's Generalized Jacobian of $\Psi$). For this reason, we turn our attention to $\mathbf{Y}$.} Thanks to Lemma \ref{lemmapseudoaffine}  and  Proposition \ref{proppseudoaffine} (with $\mathbf{Y}$ in the role of $F$) we can construct  a $QDQ$ at zero  of $\mathbf{Y}$  once   we know, for every $i=1,\ldots N$,  a $QDQ$ at zero  of its restriction $\mathbf{Y}(\varepsilon_i \mathbf{e}_i)$ to the axis $\xR\mathbf{e}_i$.  
	
	So, let us fix $i=1,\ldots N$, and let us observe that, for any $\varepsilon_i >0$ sufficiently small, one has 
	\bel{compositio}\mathbf{Y}(\varepsilon_i \mathbf{e}_i) = \Phi_{s_i}^{\ol{S}}\circ \mathbf{Y}_{s_i} ({\varepsilon_i}) \eeq
	where 
	$$ \mathbf{Y}_{s_i}({\varepsilon_i}) :=
	\Bigg(y^0_{{\varepsilon_1}}(s_i),
	y_{{\varepsilon_1}}(s_i),y^l_{{\varepsilon_1}}(s_i) \Bigg)
	%
	% \extpbf{y}(\ol{S})+\Big(0,\bm{0},\Psi\big((y^0_{\bm{{\varepsilon_1}}},y_{\bm{{\varepsilon_1}}})(\ol{S})\big)\Big).
	$$
while, for every $ q\in \xR^{1+n}$,   $[s_i,\ol{S}]\ni s\mapsto \Phi_{s_i}^{s}(q)$
	denotes   the solution to the Cauchy problem \begin{equation}\label{var2}
		\begin{cases}
			\ds\frac{d}{ds}(y^0,y,y^l)=\cl{F}\left(y,\ol{w}^0_{{\varepsilon}_i},\ol{w}_{{\varepsilon}_i}, \ol{\alpha}_{{\varepsilon}_i}\right) \cdot \Big(1+\zeta_{{\varepsilon}_i}\Big)\\
			(y^0,y,y^l)(s_i)=q,
		\end{cases}
	\end{equation}
	
	Therefore we can  again apply the chain rule for $QDQ$s (Prop. \ref{chain}) to the composed map $ \Phi_{s_i}^{\ol{S}}\circ \mathbf{Y}_{s_i} ({\varepsilon_i})$.
	
	Let us begin with determining  a $QDQ$ of $({\varepsilon_i})\mapsto\mathbf{Y}_{s_i}$ at $\varepsilon_i=0$.  We distinguish  the case when $\mathbf{c_i}$ is a needle variation generator from the one in which $\mathbf{c_i}$ is a bracket-like variation generator.
	\begin{itemize}
		\item 
		If 
		$\mathbf{c_i}=(\hat w^0,\hat w,\hat a,\hat\zeta)\in \mathfrak{V}_{ndl}$,
		standard  arguments imply  that  $(v^0_{\bm c_i,s_i},v_{\bm c_i,s_i},v^{l}_{\bm c_i, s_i})$ is the right derivative at $0$ of the path $\varepsilon_i\mapsto \mathbf{Y}_{s_i} ({\varepsilon_i})$. Therefore, the singleton $\left\{(v^0_{\bm c_i,s_i},v_{\bm c_i,s_i},v^{l}_{\bm c_i, s_i})\right\}$  a   $QDQ$ at $0$ of $\mathbf{Y}_{s_i}$ in the direction of $\xR_+$.

		\item 
		Instead, if  $\mathbf{c_i}=(j,k)\in  \mathfrak{V}_{brk}$, for some $j,k=1,\ldots,m_1$, by applying  a result established in  \cite{AngrisaniRampazzoQDQ} which proves that set-valued Lie brackets are $QDQ$s of commutator-like multiflows, we get again that the set $(v^0_{\bm c_i,s_i},v_{\bm c_i,s_i},v^{l}_{\bm c_i, s_i})$  a   $QDQ$ at $0$ of $\mathbf{Y}_{s_i}$ in the direction of $\xR_+$.
	\end{itemize}	 
	
	  Finally, by  invoking Lemma \ref{variazionale}, we get that  $$\Big\{\mathcal{E}'_k(0,0,{M},\omega), \text{ where } (M,w) \text{ is a meas. selection of } \partial_y(\ol{F}^e,\ol{l}^e)\Big\}$$ 
		is  a $QDQ$ of  $\Phi_{s_i}^{\ol{S}}$ at $\mathbf{Y}_{s_i}(0)$ in the direction of $\xR^n$.
	
	The proof of this first part of the statement of the theorem follows by use of the chain rule for Quasi-differential-quotients, namely Proposition \ref{chain}, for the reasons explained at the beginning and with the simple observation that $$\begin{pmatrix}1 & \mathbf{0} &0 \\ \mathbf{0} & \mathbf{1} & \mathbf{0} \\ m_t & m_x & 1 \end{pmatrix}\mathcal{E}'_k(0,0,{M},\omega)=\mathcal{E}'_k(m_t,m_x,{M},\omega)$$
	
	The proof of the part concerning the special case when ${\bf c}_r\in \mathfrak{V}_{ndl}$ for all $r\in 1,\ldots,N$, only requires that one extends the analysis to the last component, $\beta_{\mathbf\varepsilon}(S)$,  observing that $$\beta_{\mathbf\varepsilon}(\ol{S})-\ol{\beta}(\ol{S}):=\sum\limits_{k=1}^N \int_{s_k-{\varepsilon}_k}^{s_k} |w_k|(1+\zeta_k)-|\ol{w}(\sigma)|\,d\sigma.$$
	So the proof of the theorem is concluded in  view of the  "set product rule" in Proposition \ref{basicprops} .
\end{proof}

\subsection{Linear separability of approximating cones at the end-time}
We will use the fact that $\Lambda'_N \xR_+^N$ is a $QDQ$-approximating multicone  of  $\mathcal{R}'_\delta$ at $\big(\ol{y}^0, \ol{y}, \ol{y}^l(\ol{S}) +\ol{\Psi}(\ol{S})\big)$ to deduce  a linear separability result at time $\ol{S}$.
\begin{lmm}\label{pmend} 	Let $\opprocess$ be a canonical local minimizer for the extended  problem  $(P_{ext})$. For some positive integer  $N$, let $\Lambda'_N$ be defined as in the former subsection, and  assume that $\ol{\beta}(\ol{S})<K$ as soon as  ${\bf c}_k\in\mathfrak{V}_{brk}$ for some  $k\in\{1,\ldots,N\}$. Then:
	\begin{enumerate}
		\item for any $QDQ$-approximating multicone $\mathbcal{T}$  to the target $\mathfrak{T}$ at $\big(\ol{y}^0, \ol{y}\big)(\ol{S})$, there exist   %$(\xi_0,\xi,\xi_c)\in (\xR\times\xR^n\times\xR)^*\backslash\{(0,0,0)\}$,
		$ L' \in \Lambda'_N$, $\mathcal{T} \in \mathbcal{T}$, and $(\xi_0,\xi,\xi_c) \in \Big(L'\xR_+^N\Big)^\perp$ verifying $\xi_c\le 0$ and $(\xi_0,\xi) \in -\mathcal{T}^\perp;$
		\item
		furthermore, if $c_k \in \mathfrak{V}_{ndl}$ for all $k \in \{1,\ldots, N\}$, then the above linear form $(\xi_0,\xi,\xi_c)$  can be chosen so that
		$(\xi_0,\xi,\xi_c,\pi) \in \Big(L\xR_+^N\Big)^\perp,$  \text{ for some } $L \in \Lambda_N$ and  $\pi\leq 0$.
	\end{enumerate}
\end{lmm}
\begin{proof} By Lemma \ref{setsep} we know that, for $\delta>0$ sufficiently small,  the projected profitable set $ \mathfrak{P}^{'}$ and the projected $\delta$-reachable set  $ \mathfrak{R}^{'}_{\delta}$  are locally separated. Moreover $\Big\{\mathcal{T}\times (-\infty,0): \mathcal{T} \in \mathbcal{T}\Big\}$  is an $QDQ$-approximating multicone to the projected profitable set $ \mathfrak{P}^{'}$ at $\big(\ol{y}^0, \ol{y}, \ol{y}^l(\ol{S}) +\ol{\Psi}(\ol{S})\big)$. Hence, by Lemma \ref{OpenMappingConsequence}, it follows that the $QDQ$-approximating multicones $\Big\{\mathcal{T}\times (-\infty,0): \mathcal{T} \in \mathbcal{T}\Big\}$ and $\Lambda'_N \xR_+^N$ are not strongly transverse. Now, since $\Big\{\mathcal{T}\times (-\infty,0): \mathcal{T} \in \mathbcal{T}\Big\}$ is a multicone whose elements are contained in the semispace $\mathbb{R}^{1+n} \times \xR_-$,  by Lemma \ref{nontrasversalitaforte} we can infer the existence of 
	$(\xi_0,\xi,\xi_c)\in (\xR\times\xR^n\times\xR)^*\backslash\{(0,0,0)\}$, $ L \in \Lambda_N$, and  $\mathcal{T} \in \mathbcal{T}$,  such that $$(\xi_0,\xi,\xi_c)\in\Big(L\xR_+^N\Big)^\perp, \quad (\xi_0,\xi,\xi_c)\in -\Big(\mathcal{T}\times \xR_-\Big)^\perp. $$
	In particular, $\xi_c\leq 0$ and  $(\xi_0,\xi) \in -\mathcal{T}^\perp$, so that (1) is proved. The existence of a $\pi\le 0$ such that property (2) holds true comes from the same argument as soon as one  considers the local separation of the profitable set $ \mathfrak{P}$  and the $\delta$-reachable set $ \mathfrak{R}_{\delta}$ (in the augmented space $\xR^{1+n+1+1}$).
\end{proof}

\subsection{A maximum principle for  finitely many variations}
By using propagation due to the adjoint  differential inclusion, as a direct consequence of  Lemma \ref{pmend} we get 
a {\it maximum principle  for the instants ${s}_1,\ldots, {s}_N$ and the variation generators  $\mathbf{c}_1,\ldots, \mathbf{c}_N$:  }

\begin{lmm} 	Let $\opprocess$ be a canonical local minimizer for the extended  problem  $(P_{ext})$, and let  ${s}_1,\ldots, {s}_N\in (0,\ol{S})$,  $\mathbf{c}_1,\ldots, \mathbf{c}_N\in \mathfrak{V} $ be as above, for some integer $N>0$.
	%, with our hypotheses on $g_i,f,l,\hat{l}_1,\Psi,\mathfrak{T},A$ and $C$. 
	Let $\mathbcal{T}$ be a $QDQ$-approximating multicone for the target set $\mathfrak{T}$ at $(\ol{y}^0,\ol{y})(\ol{S})$.
	In the event that $\mathbf{c}_k \in \mathfrak{V}_{brk}$ for some $k \in \{1,\ldots, N\}$,  assume also that  $\ol{\beta}(\ol{S})<K$ and $\hat{l}_1(\cdot,0,\cdot)\equiv 0$ .
	%, and let us assume that $\ol{\beta}(\ol{S})<K$ whenever ${\bf c}_k\in\mathfrak{V}_{ndl}$ for some $k\in\{1,\ldots,N\}$
	Then, there exist $$(p_0,p,\lambda)\in \xR^*\times AC\left([0,\ol{S}], (\xR^n)^*\right)\times \xR^* \quad\,\,\, \text{ and }\,\,\,\quad\mathcal{T} \in \mathbcal{T}$$ such that $\lambda\ge 0$ and:
	\begin{itemize} 
		\item[\rm i)] {\rm \sc  (non triviality)} $(p_0,p,\lambda)\neq 0$;
		\item[\rm ii)]{\rm \sc (adjoint differential inclusion)}\bel{prop3}\frac{dp}{ds}\in -\partial_y^C H(\ol{y},p_0,p,\lambda,\pi,\ol{w}^0,\ol{w},\ol{\alpha});\eeq
		% Moreover, inequality \eqref{mp1} (equiv. (ii) of Lemma \ref{pmend} can be written as  \begin{equation}\label{prop1} p_0V^0_k+p({s}_k)  V_k-\lambda V^l_k \le 0, \quad \text{ for all } k \le N.\end{equation}
		\item[\rm iii)] {\rm \sc (non tranversality)}
		\begin{equation}\label{prop2}
			(p_0,p(\ol{S}))\in -\lambda\partial_{(t,x)}^C\Psi\Big((\ol{y}^0,\ol{y})(\ol{S})\Big) -\mathcal{T}^\perp  \,\,  ;\end{equation}
		\item[\rm iv)] {\rm \sc (first order maximization)} if 
		% ${c}_k \in \mathcal{V}_1$,
		$\mathbf{c}_k=(w^0_k,w_k,a_k,\zeta_k)\in \mathfrak{V}_{ndl}$, 
		\begin{multline}\label{prop4}
			H(\ol{y}(s_k),p_0,p(s_k),\lambda,0,w^0_k,w_k,a_k) \le \\ \le H(\ol{y}(s_k),p_0,p(s_k),\lambda,0,\ol{w}^0(s_k),\ol{w}(s_k),\ol{\alpha}(s_k))
		\end{multline}
		\item[\rm v)]{\rm \sc (nonsmooth Goh condition)} if 
		$\mathbf{c}_k=(i_k,j_k)\in \mathfrak{V}_{brk}$, 
		{ \begin{equation}\label{prop5}
			\min \,\,p(s_k) \cdot   [g_{i_k},g_{j_k}]_{set}(\ol{y}(s_k))\leq 0 \footnote{Here we mean that the minimization is performed over the elements of $$p(s_k)\cdot [g_{i_k},g_{j_k}]_{set}(\ol{y}(s_k)):=\big\{ p(s_k)V, \,\, V \in [g_{i_k},g_{j_k}]_{set}(\ol{y}(s_k))\big\}.$$}
			\end{equation}}
		
	\end{itemize}
	If, instead, $\ol{\beta}(\ol{S})=K$ and all $\mathbf{c}_k\in \mathfrak{V}_{ndl}$ for every $k=1,\ldots,N$, then there exists a triple  $(p_0,p,\lambda)$ and a real number $\pi \le 0$ such that, {\rm i)-\rm iii)}  are verified , while  inequality \eqref{prop4} is replaced by  \begin{multline}\label{prop4-bis}
		H\Big(\ol{y}(s_k),p_0,p(s_k),\lambda,\pi,w^0_k,w_k,a_k\Big) \le \\ \le H\Big(\ol{y}(s_k),p_0,p(s_k),\lambda,\pi,\ol{w}^0(s_k),\ol{w}(s_k),\ol{\alpha}(s_k)\Big)\qquad \forall k=1,\ldots,N.
	\end{multline}
\end{lmm}
\begin{proof}
	Let us observe that we can rephrase (ii) from Lemma \ref{pmend} by saying   that  there exists a linear form  $(\xi_0,\xi,\xi_c)\in (\xR\times\xR^n\times\xR)^*\backslash\{(0,0,0)\}$,  measurable selections  $M(s)\in \partial^C_y\left(f(\ol{y}(s),\ol{\alpha}(s))\ol{w}^0(s)+\sum\limits_{i=1}^m g_i(\ol{y}(s))\ol{w}^i(s)\right)$,  $\omega(s) \in \partial^C_y l^e(\ol{y}(s),\ol{w}^0(s),\ol{w}(s),\ol{\alpha}(s))$.
	 a.e. $s \in (0,\ol{S})$, and a choice of $$\left(V_j^0,V_j ,V^l_j\right)\in \left(v_{\mathbf{c}_j,{s}_j}^0,v_{\mathbf{c}_j,{s}_j},v_{\mathbf{c}_j,{s}_j}^l \right),\quad \forall  j=1,\ldots,N,$$ $$(m_t,m_x) \in \partial_{(t,x)}^C\Psi \Big(\ol{y}^0(\ol{S}),\ol{y}(\ol{S})\Big)$$such that $\xi_0\leq 0$ and, $\forall k=1,\ldots,N,$ 
	\begin{equation}\label{mp1}\begin{array}{l}\ds\xi_0V^0_k+\xi  e^{\int^{\ol{S}}_{{s}_k} M(s)}  V_k+\xi_c\Bigg[m_t V_k^0+\\\qquad\qquad\ds+m_x e^{\int_{{s}_k}^{\ol{S}} M(s)}    V_k+{\bigintsss_{{s}_k}^{\ol{S}}\omega(s)   e^{\int_{s}^{{s}_k} M(\sigma)\,d\sigma}\,ds}  V_k+V^l_k\Bigg]\le 0.\end{array}\end{equation}
	
	Setting  $\lambda:=-\xi_c$, $p_0:=\xi_0-\lambda m_t$ and, for all $s\in [0,\ol{S}]$, $$p(s):=\left(\xi-\lambda m_x\right)e^{\int_s^{\ol{S}} M(\sigma)\,d\sigma}-\lambda {\int_{s}^{\ol{S}}\omega(\sigma)   e^{\int_{\sigma}^{s} M(\tau)\,d\tau}\,d\sigma}, $$
	%{\fra lo $\xi$ qui sopra è $\xi_0$? Se sì, non sarebbe meglio
	%	$p(s):=p_0e^{\int_s^{\ol{S}} M(\sigma)\,d\sigma}-...$ ?  } {\angr no, $p_0=\xi_0-\lambda m_t$, qui abbiamo $\xi-\lambda m_x$ (ed è giusto così), che è diverso in due punti, il pedice di xi ed il pedice di m.}
	
	we get that $p(\cdot)$ satisfies the adjoint differential equation  $\ds\dot{p}(s)=-p(s)  M(s)+\lambda\omega(s).$
	%for some measurable section of the multivalued function $\partial_y\ol{F}^e$, or, in other words, 
	In particular,  $p(\cdot)$ verifies   the adjoint  differential inclusion \eqref{prop3}.
	Therefore, inequality \eqref{mp1} can be written as  \begin{equation}\label{prop1} p_0V^0_k+p(\ol{s}_k)  V_k-\lambda V^l_k \le 0, \end{equation}
	while (1) of Lemma \ref{pmend} now reads as 
	\eqref{prop2}.
	Specializing \eqref{prop1} to bracket-like variations $\mathbf{c}_k=(i_k,j_k)\in\mathfrak{V}_{brk}$ , we obtain \eqref{prop5}, whereas, when $\mathbf{c}_k=(w^0_k,w_k,a_k,\zeta_k)\in\mathfrak{V}_{ndl}$, %with $k\in \mathcal{I}_1$ 
	we get \eqref{prop4}.
	
	The case $\ol{\beta}(\ol{S})=K$ and all variations verify  $\mathbf{c}_k=(w^0_k,w_k,a_k,\zeta_k)\in\mathfrak{V}_{ndl}$ is proved similarly, by making use of (1) instead of (2) from Lemma \ref{pmend}.
	%\footnote{Actually, no proof would be needed for needle variations.}
	
\end{proof}
\subsection{Infinitely many variations} To complete  the proof  of Theorem \ref{TeoremaPrincipale}, we now combine a standard procedure,  based on Cantor's non-empty  intersection theorem,  with the crucial fact that the set-valued brackets are convex-valued.{ We will only deal  with the case when $\ol{\beta}(\ol{S})<K$,since  the case $\ol{\beta}(\ol{S})=K$ is nothing but the standard first order maximum principle applied to the rescaled, reparametrized problem.\footnote{Of course one also needs that the tangent object to the target happens to be a $QDQ$ approximating cone.} }

Begin with observing that Lusin's Theorem implies that  there exists a sequence of subsets $E_q\subset [0,\bar S]$, $q\geq 0$, such that
$E_0$ has null measure,
for every $q>0$  $E_q$ is  a compact set such that the restriction to $E_q$ of the map
$s\mapsto\left(\ol{w}^0, f(\ol{y},\ol{\alpha},\ol{w}^0)+\sum\limits_{i=1}^m g_i(\ol{y})\ol{w}^i,  l^e(\ol{y},\ol{w}^0,\ol{w}),|\ol{w}|\right)(s)$ is continuous, and
$(0,\ol{S})_{Leb}=\displaystyle\bigcup\limits_{q=0}^{+\infty} E_q.$ 
For every $q>0$ let use $D_q\subseteq E_q$  to denote  the set of all density points of $E_q$\footnote{A point $x$ is called a density point for a Lebesgue-measurable set $E$ if $\ds\lim\limits_{\rho \to 0} \frac{|B_\rho(x)\cap E|}{|B_\rho(x)|}=1$}, which,  by  Lebesgue Theorem  has the same Lebesgue measure as $E_q$. In particular, the subset  $D:=\bigcup\limits_{q=1}^{+\infty} D_q$ has measure equal to $\ol{S}$.\\
\begin{dfntn}Let $X\subseteq D\times \mathfrak{V}$ be any subset of time-generator pairs. We will say that a triple $(p_0,p,\lambda)\in \xR \times AC([0,\ol{S}];\xR^n)\times \xR_+$ {\it satisfies property $(P_X)$} if the following conditions  (1)-(3) are verified:
	\begin{itemize}
		\item[(1)] $p$ is a solution on $[0,\ol{S}]$ of  the differential inclusion \begin{equation} \label{differentialinclusion}\dot{p}\in -p \,\, \partial_y^C\left(f(\ol{y},\ol{\alpha})\ol{w}^0+\sum\limits_{i=1}^m g_i(\ol{y})\ol{w}^i\right)+\lambda\partial_y^C l^e(\ol{y},\ol{w}^0,\ol{w},\ol{\alpha});\end{equation} 
		\item[(2)] one has  \begin{equation}\label{propnontrasversalita} (p_0,p(\ol{S}))\in -\lambda \partial_{(t,x)}^C\Psi \Big(\ol{y}^0(\ol{S}),\ol{y}(\ol{S})\Big)  -\ol{\bigcup\limits_{\mathcal{T}\in \mathbcal{T}} \mathcal{T}^\perp}   ;\end{equation}
		\item[(3)] for every $(s,\mathbf{c}) \in X$, if $\mathbf{c}=(w^0,w,a,\zeta)\in\mathfrak{V}_{ndl}$, then \begin{multline} \label{disequazione} p_0w^0(1+\zeta)+p(s)\,\left(f(\ol{y}(s),a){w}^0+\sum\limits_{i=1}^m g_i(\ol{y}(s)){w}^i\right)(1+\zeta)-\lambda l^e(\ol{y}(s),w^0,w,a)\le \\  p_0\ol{w}^0
			%(1+\ol{\zeta})
			+p(s)\,\left(f(\ol{y}(s),\ol{\alpha}(s))\ol{w}^0(s)+\sum\limits_{i=1}^m g_i(\ol{y}(s))\ol{w}^i(s)\right)
			%\left(1+\ol{\zeta}(s)\right)
			-\\
			\qquad\qquad\qquad\qquad\qquad\qquad\qquad\qquad-\lambda l^e(\ol{y}(s),\ol{w}^0(s),\ol{w},\ol{\alpha}(s)), \end{multline} while, if $\mathbf{c}=(i,j)\in \mathfrak{V}_{brk}$, then \begin{equation}\label{condizionealsecondoordine}\displaystyle\min\limits_{\qquad V \in [g_i,g_j]_{set}(\ol{y}(s))} p_n(s)  V\le 0.\end{equation}  
	\end{itemize}

	Finally, for any given  $X \subseteq D\times \mathfrak{V}$, let us  define the subset $\Theta(X)\subset \xR^*\times AC\big([0,\ol{S}];(\xR^n)^*\big)\times \xR^*$ as $$\Theta(X):=\left\{\begin{aligned} & (p_0,p,\lambda) \in \xR\times AC([0,\ol{S}];\xR^n)\times \xR\,:\, |(p_0,p(\ol{S}),\lambda)|=1,\,\\ & (p_0,p,\lambda) \text{ verifies the property } (P_X)\end{aligned}\right\}.$$
\end{dfntn}\begin{lmm}{\it For any subset  $X \subseteq D\times \mathfrak{V}$, $\Theta(X)$ is a compact  subset of $\xR\times AC\big([0,\ol{S}];\xR^n\big)\times \xR$}, when the latter  is endowed  with the norm $\|(p_0,p(\cdot),\lambda)\|:=|p_0|+\|p\|_\infty+|\lambda|$. \end{lmm}
\begin{proof} Consider a sequence $(p_{0,n},p_n(s),\lambda_n) \in \Theta(X)$. The set-valued maps
	$$ s\mapsto \partial_y^C\left(f(\ol{y}(s),\ol{\alpha}(s))\ol{w}^0(s)+\sum\limits_{i=1}^m g_i(\ol{y}(s))\ol{w}^i(s)\right)
	$$
	$$\ds s\mapsto \partial_y^C l^e\Big(\ol{y}(s),\ol{w}^0(s),\ol{w}(s),\ol{\alpha}(s)\Big)$$
	have uniformly bounded closed convex values as they are Clarke Jacobians of functions that are globally Lipschitz (after the non-restrictive  cut off operation described earlier).
	% to a compact set containing a tubular neighbourhood of the set $\ol{y}([0,\ol{S}])$.
	%
	% For every $s\in [0,\ol{S}]$, $\partial_y\left(f(\ol{y}(s),\ol{\alpha}(s))\ol{w}^0(s)+\sum\limits_{i=1}^m g_i(\ol{y}(s))\ol{w}^i(s)\right)$ $\frac{\partial l^e}{\partial y}(\ol{y}(s),\ol{w}^0(s),\ol{w}(s)) $  is a non-empty convex compact set of vectors contained in the ball of radius the  Lipschitz constant of $F^e$, while $\ds\frac{d\ol{l}^e}{dy}$  is bounded by the maximum of the continuous function $\ds\frac{dl^e}{dy}$ over some compact set containing a tubular neighbourhood of the reference trajectory.
	Furthermore, the quantities $|p_n(\ol{S})|$, $\lambda_n$ and $p_{0,n}$ are bounded in norm by $1$, so that  we are in the position to use the following fact (which can  be deduced, e.g., from Theorem 1 in Chapter 2 of \cite{Aubin}):
	
	\begin{itemize}\item{\it 
			Let $C(s):[0,\ol{S}]\setmap \xR^n$ and $B(s):[0,\ol{S}]\setmap \xR^n$ be a measurable set-valued map with compact, convex, non-empty values. Moreover, assume that there exists $R>0$ such that, for every $s\in [0,\ol{S}]$, the sets $B(s)\subset Lin(\xR^n,\xR^n)$ are $C(s)\subset\xR^n$ are all contained in the  ball centered in the  corresponding origins  and of radius $R$.
			Let $p_n(s)$ be a sequence of solutions to the differential inclusion
			\begin{equation}\label{inclusionedifferenzialstantard} \dot{p}(s)\in p(s)B(s)+C(s),\qquad \text{ for almost all } s \in [0,\ol{S}]\end{equation} all satisfying $|p_n(\ol{S})|\le 1$. Then there is a subsequence of $p_n(s)$ that uniformly converges to a function $p(s)$, and $p(s)$ is also a solution to the differential inclusion \eqref{inclusionedifferenzialstantard}.
		}
	\end{itemize}

	Therefore,  modulo thrice extracting subsequences from our sequence, we can assume \begin{center}$\lambda_n\to \lambda \ge 0$, $p_{0,n}\to p_0$ and $p_n\to p\in AC$, uniformly for $s \in [0,\ol{S}]$,\end{center} with $p(s)$ still satisfying the differential inclusion \eqref{differentialinclusion}.  Since the paths  $p_n$ converges uniformly to $p$, properties  \eqref{disequazione} and \eqref{propnontrasversalita} are  inherited by $p(s)$ from the sequence $p_n(s)$ by passing to the limit. Finally,   passing to the limit we get that  \eqref{condizionealsecondoordine} holds true as well.\end{proof}
%, as it is known that $$x\to \min_{y \in M} f(x,y)$$ is a continuous function whenever $$f:(x,y)\in H\times M \mapsto f(x,y)\in \xR$$ is a continuous function on a compact set $H\times K$

\begin{lmm}\label{lemnone}{The set $\Theta(D\times \mathfrak{V})$ is non-empty}.\end{lmm}
\begin{proof} We are going to use a non-empty intersection argument, which is quite standard, except for  the part concerning the  set-valued bracket. Let us notice that $$\Theta(X_1\cup X_2)=\Theta(X_1)\cap \Theta(X_2),
	\quad \forall X_1,X_2\subseteq D \times \mathfrak{V},$$ so that
	\begin{equation}\label{lemint}\Theta(D\times \mathfrak{V})=\bigcap\limits_{\substack{X \subseteq D \times \mathfrak{V} \\X \text{ finite }}} \Theta(X).\end{equation}
	We have to prove that this infinite intersection is non-empty. 
	We begin with proving the following fact: 
	
	{\bf Claim 1} {The set \it $\Theta(X)$ is non-empty as soon as  $X$ is finite.} 
	\\\,
	Indeed, by \eqref{prop1}-\eqref{prop5} we already know that  that $\Theta(X)\neq\emptyset$  whenever $X$ comprises $N$ couples $(s_k,\mathbf{c}_k)\in D\times \mathfrak{V}$ such that $s_k<s_l$ whenever $1\le k<l\le N$.
	We have to show  that we can allow  $X$ to have the general  form 
	$$X=\Big\{ (s_k,\mathbf{c}_k)\in D\times \mathfrak{V},  \quad s_k\leq s_l\,\,\text{as soon as}\,\, 1\le k<l\le N\Big\}.$$  For any $(r,k)\in\mathbb{N}\times\{1,\ldots,N\}$, choose an instant   $s_{k,r}\in E:=\bigcup\limits_{q=1}^\infty E_q$ in such a way that such that $s_{1,r}<\ldots<s_{N,r}$ and the sequences  $(s_{k,r})_{r\in\mathbb{N}}$, $1\leq k\leq N$,  converge to $s_k$. For every $r\in \mathbb{N}$, Consider  the sets $X_r:=\{(s_{k,r},\mathbf{c}_k),\,\,k\le N\}$ , $r\in\mathbb{N}$.
	By the previous steps we know that  $\Theta(X_r)$ is non-empty, so that we can choose  $(p_{0,r},p_r(s),\lambda_r) \in \Theta(X_r)$.
	% We want to show now that, as $r\to +\infty$, they converge to $(p_0,p(s),\lambda) \in \Theta(X)$. 
	Again, modulo thrice extracting subsequences, we can assume that \begin{enumerate}\item[i] $p_{0,r}$ converges to a real number $p_0$, \item[ii] $p_r$ uniformly converges to an absolutely continuous function $p$ solving the adjoint differential inclusion, and \item[iii] $\lambda_r$ converges to a non-negative real number $\lambda$.\end{enumerate} Therefore,   $(p_0,p)$  inherits the non-transversality condition \eqref{propnontrasversalita}. Moreover, the  uniform convergence of the $p_r$ and the continuity of the involved functions  imply that $(p_0,p,\lambda)$  verifies the {  Hamiltonian maximization} \eqref{disequazione}  at the instants $s_k$.\\
	What is going to be a little trickier to prove is that, if $\mathbf{c}_k=(i_k,j_k)$ is a bracket-like variation generator, \eqref{condizionealsecondoordine} holds at time $s_k$ for the multiplier $p$, starting from the fact that  it is satisfied at time $s_{k,r}$ by the multiplier $p_r$. As a matter of fact, what we know is that, for every natural number $r$,  $$p_r(s_{k,r}) V_{k,r} \le 0$$ for some $V_{k,r}  \in [g_{i_k},g_{j_k}]_{set}(\ol{y}(s_{k,r})).$ By definition of the set-valued Lie-bracket as closure of a convex hull, this implies that, for every $r\in\mathbb{N}$, there is a sequence $(y_{k,r,n})_{n\in\mathbb{N}}$ of differentiability points for both $g_{i_k}$ and $g_{j_k}$ converging to $\ol{y}(s_{k,r})$ and  such that $$p_r(s_{k,r})  \left(\lim\limits_{n\to +\infty} Dg_{i_k}(y_{k,r,n})g_{j_k}(y_{k,r,n})-Dg_{j_k}(y_{k,r,n})g_{i_k}(y_{k,r,n})\right) \le \frac{1}{r}.$$ Therefore, there is a large enough number $N_{1,r}$,  such that, for all $n\ge N_{1,r}$, $$p_r(s_{k,r})  \Big( Dg_{i_k}(y_{k,r,n})g_{j_k}(y_{k,r,n})-Dg_{j_k}(y_{k,r,n})g_{i_k}(y_{k,r,n})\Big) \le \frac{2}{r}.$$ Also, since $\lim\limits_{n\to +\infty} Dg_{i_k}(y_{k,r,n})g_{j_k}(y_{k,r,n})-Dg_{j_k}(y_{k,r,n})g_{i_k}(y_{k,r,n})$ is bounded  as $r$ varies in $\mathbb{N}$, we can assume, modulo extracting a subsequence, that the limit \begin{equation}\label{limitex}\lim\limits_{r\to +\infty}\lim\limits_{n\to +\infty} Dg_{i_k}(y_{k,r,n})g_{j_k}(y_{k,r,n})-Dg_{j_k}(y_{k,r,n})g_{i_k}(y_{k,r,n})
	\end{equation} does exist. Let us call $V_k$ this limit, i.e. let us set 
	$$V_k:=\lim\limits_{r\to +\infty}\lim\limits_{n\to +\infty} \big[g_{j_k}, g_{i_k}\big](y_{k,r,n}).$$
	Since $\ds\lim_{n\to\infty}y_{k,r,n}= \ol{y}(s_{k,r})$ for any $r$    and  $\ds  \lim_{r\to\infty}\ol{y}(s_{k,r}) =\ol{y}(s_k)$, we can construct a sequence  $\ds (N_{2,r})_{r\in\mathbb{N}}$ of natural numbers  such that \begin{enumerate}\item[i]$N_{2,r}\to +\infty $, \item[ii] $y_{k,r,N_{2,r}}\to \ol{y}(s_k)$, and \item[iii] $y_{k,r,N_{2,r}}$ is a  point of differentiability for both $g_{i_k}$ and $g_{j_k}$ for any $r\in\mathbb{N}$.\end{enumerate} By the existence of the limit \eqref{limitex} we deduce that,  by taking, for any $r$,  a suitably large $N_{3,r}>N_{2,r}>N_{1,r}$  one has 
		$$V_k:=\lim\limits_{r\to +\infty}\lim\limits_{n\to +\infty} \big[g_{j_k}, g_{i_k}\big](y_{k,r,N_{3,r}}).$$
	As $N_{3,r}>N_{2,r}$, this means $V_k \in [g_{i_k},g_{j_k}]_{set}(\ol{y}(s_k))$, in that  it is the limit of the Lie bracket of $g_{i_k}$ and $g_{j_k}$ computed along a sequence of points that converges to $\ol{y}(s_k)$. Moreover, since $N_{3,r}>N_{1,r}$,  one has
		$$p_r(s_{k,r})\cdot \big[g_{j_k}, g_{i_k}\big](y_{k,r,N_{3,r}}) \le \frac{2}{r} .$$  By passing to the limit as $r$ goes to infinity, we get  that $$p(s_k) \, V_k\le 0\qquad \Big(\text{ and } V_k \in [g_{i_k},g_{j_k}]_{set}(\ol{y}(s_k))\Big),$$  which concludes the proof of  {\bf Claim 1}.
	
	In view of  {\bf Claim 1}  and of 
	\eqref{lemint}, by Cantor's intersection theorem  we can conclude that  $\Theta(D\times \mathfrak{V})\neq\emptyset$,  for  it is the intersection of a family of compact non-empty sets with the property that any finite intersection of these sets is non-empty. Now, any $(p_0,p,\lambda) \in \Theta(D\times \mathfrak{V})$ clearly satisfies  conditions {\bf i)} -{\bf v)} of Theorem \ref{TeoremaPrincipale}. As for the  second order condition {\bf vi)}, if \eqref{condizionealsecondoordine} holds for $\mathbf{c}=(i,j)$ and $\mathbf{c}=(j,i)$ at the same time, then we know that there are $V^-, V^+ \in [g_i,g_j]_{set}(\ol{y}(s))$ such that $p(s)  V^-\le 0$ and $p(s)  V^+\ge 0.$ Since $[g_i,g_j]_{set}(\ol{y}(s))$  is convex, this implies that there exists $\alpha\in [0,1]$ such that  $V :=\alpha V^- +(1-\alpha) V^+$ (belongs to $[g_i,g_j]_{set}(\ol{y}(s))$ and) satisfies $p(s)   V=0.$
\end{proof}

\section{An example}
In  the minimum problem \eqref{ex} below, the fact that a certain control $\bar u$ is not optimal is not deduced by the standard  first order maximum principle. In other words conditions $i)$ through $iv)$ turn out to be satisfied while, instead, the nonsmooth Goh condition condition $v)$ is not verified, so  the optimality of $\bar u$ is ruled out.

Given the  function $
\Psi(t,x) := |x|^2 +(t-1)^2
$, we will consider the following Mayer optimal control problem:
%the$$$$on the family of feasible processes of the control systemin  $\xR^3$ (plus  the usual total variation variable $\nu$), 
\bel{ex} \left.\begin{array}{c}
	\min \Psi(T,x(T))\\ \text{in the set of solutions to}\\
	\left\{\begin{array}{l}\ds
		\frac{dx}{dt} = f(x) + g_1(x)u^1 + g_2(x)u^2
		\\\ds (x^1,x^2,x^3,\nu)(0)= \left(\begin{matrix}1, 0, 2 ,0\end{matrix}\right)\qquad (x^1,x^2,x^3)\in\mathfrak{T}\times [0,K],\,\,\|\nu\|_1\leq 4 \end{array}\right.\end{array}\right.
\eeq
where the control $u$ take values in $C:=\xR^2$, the target is defined as  $\mathfrak{T}:=[0,1]\times B_{\frac{1}{2}}\left(0,0,\frac{1}{2}\right)$, and
$$
f(x) :\equiv \left(\begin{matrix}0\\ 0\\ -1\end{matrix}\right),\qquad 
g_1(x) := \left(\begin{matrix}1\\ 0\\ \displaystyle  x^2 -{|x^2|}\end{matrix}\right), \qquad 
g_2(x) := \left(\begin{matrix}0\\ 1\\ \displaystyle  x^1 +{|x^1|}\end{matrix}\right).
$$
We will show that the control $\ol{u}(s)\equiv (-1,0)$ satisfies conditions $i)$ through $iv)$ but not $v)$.\\
% In contrast, at the end of the section, a global minimizer is presented.\\
The extended problem reads
$$\min \Psi(y^0(S),y(S))$$
$$\left\{\begin{array}{l}
	\left(\begin{matrix}\ds\frac{dy^0}{^{}ds} , \frac{dy^1}{ds} ,  \frac{dy^2}{ds} , \frac{dy^3}{ds} , \frac{d\beta}{ds} \end{matrix}\right)=\left(\begin{matrix} w^0 \,, \,  w^1 \, , \, w^2  \,,  \,\left(y^2 -{|y^2|}\right)w^1 + \left( y^1 +{|y^1|}\right)w^2 - w^0\,,\, |w| \end{matrix}\right) \\\\
	
	\ds 
	(y^0,y^1,y^2,y^3,\beta)(0) = \left(\begin{matrix}0, 1, 0, 2,  0\end{matrix}\right) \qquad (y^0,y^1,y^2,y^3,\beta)(\ol{S})\in\mathfrak{T}\times [0,4]
\end{array}\right.
$$

In the notation of the extended problem, we shall focus on the (constant) space-time control $(\bar w^0,\bar w^1, \bar w^2)$ defined by $\ol{S}=2$ and 
$$
\left(\bar w^0 , \bar w ^1,\bar w^2 \right) (s) := \left(\frac{1}{2} , -\frac{1}{2} ,  0\right)\quad\forall s\in [0,2] .\quad\footnote{Actually, since $w^0\equiv \frac{1}{2}>0$ all of these extended sense controls correspond to strict sense controls  $(1,-1,0)$. This will not be true for the global minimizer presented at the end of the section.}
$$
and the corresponding trajectory
$$
\left(\bar y^0 , \bar y^1,\bar y^2,\bar y^3\right) (s):=\left(\frac{1}{2}s , 1-\frac{1}{2}s ,  0 , 2-\frac{1}{2}s\right)
$$
ending at $\left(\bar y^0 , \bar y^1,\bar y^2,\bar y^3\right) (2)=\big(1,0,0,1\big)$, with final cost $\Psi\Big(\left(\bar y^0 , \bar y\right) (2)\Big)=1.$\\
Our adjoint differential inclusion reads $$\frac{dp_1}{ds} = 0, \,\, \frac{dp_2}{ds}\in [0,p_3(s)],\,\,\frac{dp_3}{ds}=0,$$ so that,  $p_1$ and $p_3$ need to be constant.
As a $QDQ$-approximating multicone to the target set, we consider the (Boltyanski) cone $\xR_- \times \xR \times \xR \times \xR_{-}$. Moreover $\Psi$ is a smooth function and  $\partial_{(t,x)}^C\Psi\Big(\left(\bar y^0 , \bar y^1,\bar y^2,\bar y^3\right) (2)\Big) = \{(0,0,0,2)\}$, so that transversality condition reads $$(p_0,p_1,p_2,p_3)(2)\in ]-\infty,0] \times \{0\} \times \{0\} \times ]-\infty,-2\lambda].$$
	Therefore $p_1\equiv 0$. Furthermore, if  $p_3$ were equal to  $0$, then $\lambda$ would also be $0$ and $p_2$ would be constantly equal to $0$. Finally, from the first order maximization condition, one would imply $p_0=0$, and the non-triviality condition would be violated. Hence $p_3<0$, so that  we can assume $p_3=-1$.\\
The transversality condition and adjoint differential inclusion together imply that, if some multipliers $(p_0,p,\lambda,\pi)$ exist satisfying conditions $i)$ through $iv)$, then $p_1 \equiv 0$, $p_3 \equiv -1$ and $p_2(s)$ is a non-increasing function ending at $0$ (more precisely, $\ds\frac{\partial p_2}{\partial s} \in [-1,0]$).\\
Writing the first order maximization condition, we get that $$H=p_0w^0+p_2(s)w^2+(-1)[(2-s)w^2-w^0]=(p_0+1)w^0+w^2[p_2(s)-(2-s)]$$ should be maximized by the choice $w^2=0$ and $w^0=\frac{1}{2}$ for almost all $s$. This is  possible only if $p_0=-1$ and $p_2(s)=2-s$ for almost all $s$ in $[0,2]$, which means $p_2(s)=2-s$ for all $s$ because $p_2$ is absolutely continuous. In other words, the adjoint path  $$\Big(\left(\bar y^0 , \bar y^1,\bar y^2,\bar y^3\right)\,,\, (p_0,p_1,p_2,p_3)\Big)(s):=\left(\left(\frac{1}{2}s , 1-\frac{1}{2}s ,  0 , 2-\frac{1}{2}s\right)\,,\,(-1,0,2-s,-1)\right)$$  satisfies the first order conditions of the maximum principle.

However, {\it the nonsmooth Goh condition is not verified}, because the set-valued Lie bracket $[g_1,g_2]_{set}(y^1,y^2,y^3)$ is equal to  $\{0\}\times \{0\} \times [2,4]$ whenever $y^1>0$ and $y^2=0$, which yields $$0 \not \in [-4,-2]=(p_1,p_2,p_3)(s)\cdot \Big(\{0\}\times \{0\} \times [2,4]\Big).$$\\
On the other hand, let us consider the (impulsive) control 
$$(\hat{w}^0,\hat{w}^1,\hat{w}^2):=\begin{cases}
	\left(\frac{1}{2},-\frac{1}{2},0\right) & \text{ if } s \in [0,2]\\
	\left(0,0,1\right) & \text{ if } s \in \,]2,2+\frac{\sqrt{2}}{2}]\\
	\left(0,1,0\right) & \text{ if } s \in \,]2+\frac{\sqrt{2}}{2},2+\sqrt{2}]\\
	\left(0,0,-1\right) & \text{ if } s \in \,]2+\sqrt{2},2+3\frac{\sqrt{2}}{2}]\\
	\left(0,-1,0\right) & \text{ if } s \in \,]2+3\frac{\sqrt{2}}{2},2+2\sqrt{2}]\\
\end{cases}$$
It is easy to see that the corresponding trajectory $(\hat{y}^0,\hat{y}^1,\hat{y}^2,\hat{y}^3)$ ends at $(1,0,0,0)$ which is a point in the target set, but also a point of global minimum of $\Psi$. The process turns out to be feasible because one also has  $\beta(\ol{S})=\beta(2+2\sqrt{2})=1+2\sqrt{2}<K=4$. 
Let us   underline a crucial difference between the geometrical pictures of the two extremals: in the latter,  it is simple to check that the $QDQ$-approximating cone $\xR_- \times \xR \times \xR \times \xR_{+}$ to the target at the end-point  $(1,0,0,0)$ allows for the adjoint map  $$(p_0,p_1,p_2,p_3)\equiv 0, \quad \lambda=1, \quad \pi=0 ;$$on the contrary, in the case of   the control $(\ol{w}^0 ,\ol{w}^1,\ol{w}^2)$, the $QDQ$-approximating cone  $\xR_- \times \xR \times \xR \times \xR_{-}$ to the target at the end-point  $(1,0,0,1)$  forced $p_3$ to be strictly negative, which contrasted with the nonsmooth Goh condition.

\end{document}